\documentclass[final]{siamltex}
\usepackage{graphicx}
\usepackage{psfrag}
\usepackage{amsmath}
\usepackage{amssymb}
\usepackage{amsfonts}
\usepackage{stmaryrd}
\usepackage{float}
\usepackage{subfig}
\usepackage{hyperref}
\usepackage{diagbox}
\usepackage{makecell}
\usepackage{multirow}
\usepackage{afterpage}
\usepackage{placeins}

\usepackage[usenames,dvipsnames]{color}
\newcommand{\FK}[1]{{#1}}
\newcommand{\MG}[1]{{#1}}

\newtheorem{assumption}{Assumption}

\begin{document}

\title{Nonlinear Preconditioning: How to use a Nonlinear Schwarz
    Method to Precondition Newton's Method}

\author{
V. Dolean\footnotemark[1] \and
M.J. Gander\footnotemark[2] \and
W. Kheriji\footnotemark[3] \and
 F. Kwok\footnotemark[4] \and
R. Masson\footnotemark[5]
}

\maketitle

\begin{abstract}
  For linear problems, domain decomposition methods can be used
  directly as iterative solvers, but also as preconditioners for
  Krylov methods.  In practice, Krylov acceleration is almost always
  used, since the Krylov method finds a much better residual
  polynomial than the stationary iteration, and thus converges much
  faster. We show in this paper that also for non-linear problems,
  domain decomposition methods can either be used directly as
  iterative solvers, or one can use them as preconditioners for
  Newton's method. For the concrete case of the parallel Schwarz
  method, we show that we obtain a preconditioner we call RASPEN
  (Restricted Additive Schwarz Preconditioned Exact Newton) which is
  similar to ASPIN (Additive Schwarz Preconditioned Inexact Newton),
  but with all components directly defined by the iterative
  method. This has the advantage that RASPEN already converges when
  used as an iterative solver, in contrast to ASPIN, and we thus get a
  substantially better preconditioner for Newton's method. The
  iterative construction also allows us to naturally define a coarse
  correction using the multigrid full approximation scheme, which
  leads to a convergent two level non-linear iterative domain
  decomposition method and a two level RASPEN non-linear
  preconditioner. We illustrate our findings with numerical results on
  the Forchheimer equation and a non-linear diffusion problem.
\end{abstract}

\pagestyle{myheadings}
\thispagestyle{plain}
\renewcommand{\thefootnote}{\fnsymbol{footnote}}
\footnotetext[1]{Department of Maths and Stats, University of Strathclyde,
  Glasgow, United Kingdom {\tt Victorita.Dolean@strath.ac.uk}}
\footnotetext[2]{Section de Math\'ematiques, Universit\'e de
Gen\`eve, CP 64, 1211 Gen\`eve, Switzerland, {\tt Martin.Gander@math.unige.ch}}
\footnotetext[3]{Laboratoire J.-A. Dieudonn\'e,
Universit\'e Nice Sophia Antipolis and project Coffee Inria Sophia Antipolis M\'editerran\'ee,
France, {\tt kheriji.walid@gmail.com}}
\footnotetext[4]{Department of Mathematics, Hong Kong Baptist
   University, Kowloon Tong, Hong Kong, {\tt felix$\_$kwok@hkbu.edu.hk}}
\footnotetext[5]{Laboratoire J.-A. Dieudonn\'e,
Universit\'e Nice Sophia Antipolis and project Coffee Inria Sophia Antipolis M\'editerran\'ee,
France, {\tt Roland.Masson@unice.fr}}
\section{Introduction}

Non-linear partial differential equations are usually solved after
discretization by Newton's method or variants thereof. While Newton's
method converges well from an initial guess close to the solution, its
convergence behaviour can be erratic and the method can lose all its
effectiveness if the initial guess is too far from the
solution. Instead of using Newton, one can use a domain decomposition
iteration, applied directly to the non-linear partial differential
equations, and one obtains then much smaller subdomain problems, which
are often easier to solve by Newton's method than the global
problem. The first analysis of an extension of the classical
alternating Schwarz method to non-linear monotone problems can be
found in \cite{Lions:1988:SAM}, where a convergence proof is given at
the continuous level for a minimization formulation of the problem.  A
two-level parallel additive Schwarz method for non-linear problems was
proposed and analyzed in \cite{dryja1997nonlinear}, where the authors
prove that the non-linear iteration converges locally at the same rate
as the linear iteration applied to the linearized equations about the
fixed point, and also a global convergence result is given in the case
of a minimization formulation under certain conditions. In
\cite{lui1999schwarz}, the classical alternating Schwarz method is
studied at the continuous level, when applied to a Poisson equation
whose right hand side can depend non-linearly on the function and its
gradient. The analysis is based on fixed point arguments; in addition, the
author also analyzes linearized variants of the iteration in which the
non-linear terms are relaxed to the previous iteration. A
continuation of this study can be found in \cite{lui2002linear}, where
techniques of super- and sub-solutions are used. Results for more
general subspace decomposition methods for linear and non-linear
problems can be found in \cite{xu1996two,tai1998rate}. More recently,
there have also been studies of so-called Schwarz waveform relaxation methods
applied directly to non-linear problems: see
\cite{Gander:1998:WRA,Gander:2005:OSW,descombes2011schwarz}, where
also the techniques of super- and sub-solutions are used to analyze
convergence, and \cite{halpern2009nonlinear,caetano2010schwarz} for
optimized variants.

Another way of using domain decomposition methods to solve non-linear
problems is to apply them within the Newton iteration in order to solve the
linearized problems in parallel. This leads to the
Newton-Krylov-Schwarz methods \cite{cai1994newton,cai1998parallel},
see also \cite{cai1994domain}. We are however interested in a
different way of using Newton's method here. For linear problems,
subdomain iterations are usually not used by themselves; instead, the equation
at the fixed point is solved by a Krylov method, which greatly reduces
the number of iterations needed for convergence. This can also be done
for non-linear problems: suppose we want to solve $F(u)=0$ using the
fixed point iteration $u^{n+1}={\cal G}(u^n)$. To accelerate convergence, we
can use Newton's method to solve ${\cal F}(u):={\cal G}(u)-u=0$ instead. We
first show in Section \ref{SimpleSec} how this can be done for a
classical parallel Schwarz method applied to a non-linear partial
differential equation, both with and without coarse grid, which leads
to a non-linear preconditioner we call RASPEN. With our approach, one can obtain
in a systematic fashion nonlinear preconditioners for Newton's method
from any domain decomposition method. A different non-linear
preconditioner called ASPIN was invented about a decade ago in
\cite{cai2002nonlinearly}, see also the earlier conference publication
\cite{cai2001nonlinear}. Here, the authors did not think of an
iterative method, but directly tried to design a non-linear two level
preconditioner for Newton's method.  This is in the same spirit as
some domain decomposition methods for linear problems that were
directly designed to be a preconditioner; the most famous example is
the additive Schwarz preconditioner \cite{Dryja:1987:AVS}, which does
not lead to a convergent stationary iterative method without a
relaxation parameter, but is very suitable as a preconditioner, see
\cite{gander2008schwarz} for a detailed discussion. It is however
difficult to design all components of such a preconditioner, in
particular also the coarse correction, without the help of an
iterative method in the background. We discuss in Section
\ref{ComparisonSec} the various differences between ASPIN and RASPEN.
Our comparison shows three main advantages of RASPEN: first, the one-level
preconditioner came from a convergent underlying iterative
method, while ASPIN is not convergent when used as an iterative solver
without relaxation; thus, we have the same advantage as in the linear case,
see \cite{Efstathiou:2003:WRA,gander2008schwarz}. Second, the coarse
grid correction in RASPEN is based on the full approximation scheme
(FAS), whereas in ASPIN, a different, ad hoc construction based on a
precomputed coarse solution is used, which is only good close to the
fixed point. And finally, we show that the underlying iterative method
in RASPEN already provides the components needed to use the exact
Jacobian, instead of an approximate one in ASPIN. These three
advantages, all due to the fact that RASPEN is based on a convergent
non-linear domain decomposition iteration, lead to substantially lower
iteration numbers when RASPEN is used as a preconditioner for Newton's
method compared to ASPIN.  We illustrate our results in Section
\ref{NumSec} with an extensive numerical study of these methods for
the Forchheimer equation and a non-linear diffusion problem.

\section{Main Ideas for a Simple Problem}\label{SimpleSec}

To explain the main ideas, we start with a one dimensional non-linear
model problem
\begin{equation}\label{1dmodel}
  \begin{array}{rcll}
    {\cal L}(u)&=&f,\quad&\mbox{in $\Omega:=(0,L)$},\\
    u(0)&=&0, \\
    u(L)&=&0,
  \end{array}
\end{equation}
where for example ${\cal L}(u)=-\partial_x((1+u^2)\partial_x u)$.
One can apply a classical parallel Schwarz method to solve such
problems. Using for example the two subdomains $\Omega_1:=(0,\beta)$
and $\Omega_2:=(\alpha,L)$, $\alpha < \beta$, the classical parallel Schwarz method is
\begin{equation}\label{ClassicalSchwarz1dmodel}
  \begin{array}{rcll}
    {\cal L}(u_1^n)&=&f,\quad&\mbox{in $\Omega_1:=(0,\beta)$},\\
    u_1^n(0)&=&0, \\
    u_1^n(\beta)&=&u_2^{n-1}(\beta),\\
    {\cal L}(u_2^n)&=&f,\quad&\mbox{in $\Omega_2:=(\alpha,L)$},\\
    u_2^n(\alpha)&=&u_1^{n-1}(\alpha), \\
    u_2^n(L)&=&0.
  \end{array}
\end{equation}
This method only gives a sequence of approximate
solutions per subdomain, and it is convenient to introduce a global
approximate solution, which can be done by glueing the approximate
solutions together. A simple way to do so is to select values from one of
the subdomain solutions by resorting to a non-overlapping decomposition,
\begin{equation}
  u^n(x):=\left\{\begin{array}{ll}
    u_1^n(x)&\mbox{if $0\le x <\frac{\alpha+\beta}{2}$},\\
    u_2^n(x)&\mbox{if $\frac{\alpha+\beta}{2}\le x \le L$},
  \end{array}\right.
\end{equation}
which induces two extension operators ${\widetilde P}_i$ (often called
$\tilde{R}_i^T$ in the context of RAS); we can write
$u^n={\widetilde P}_1u_1^n+{\widetilde P}_2u_2^n$.

\begin{figure}
  \centering
  \includegraphics[width=0.46\textwidth]{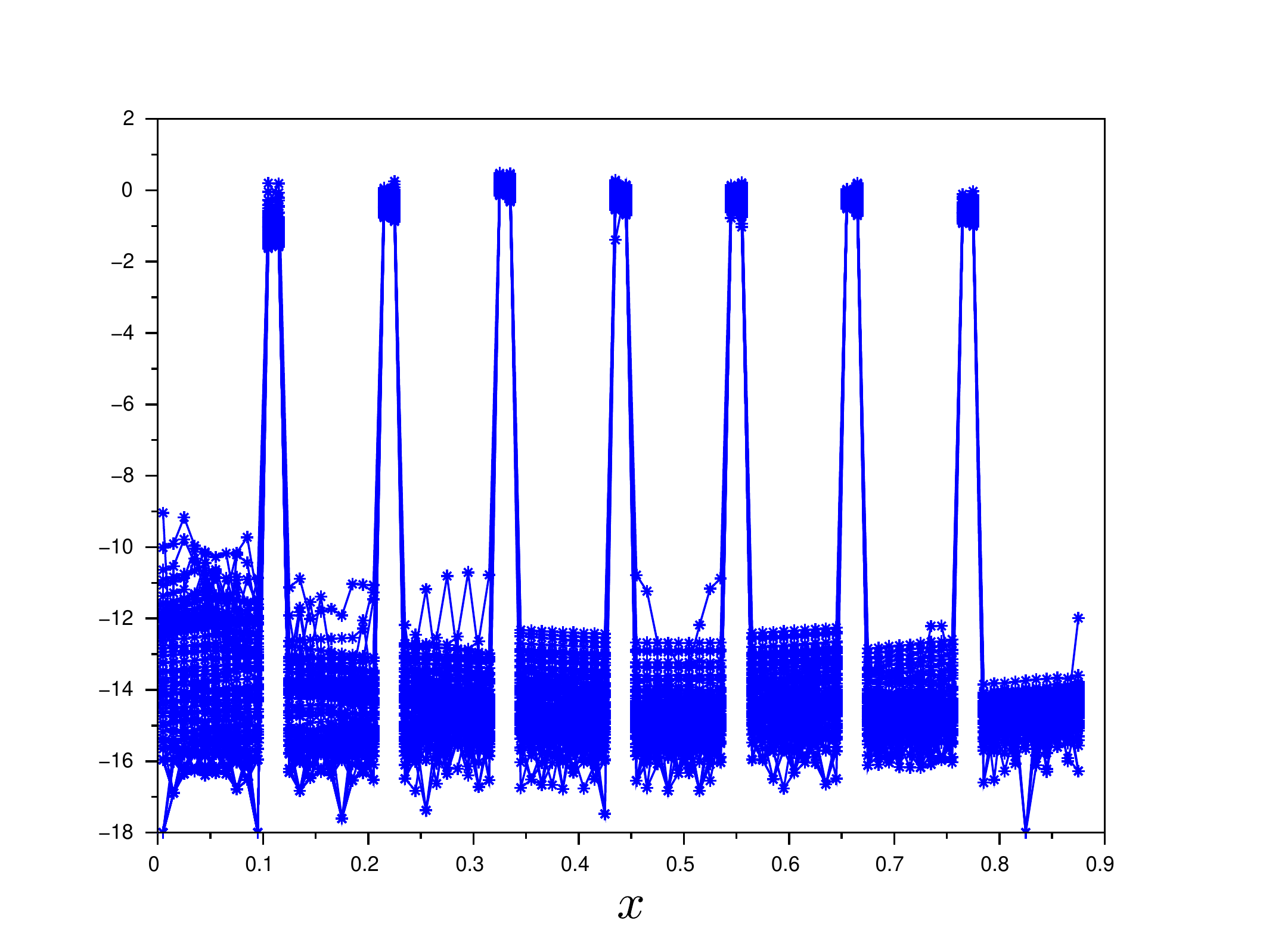}
  \includegraphics[width=0.46\textwidth]{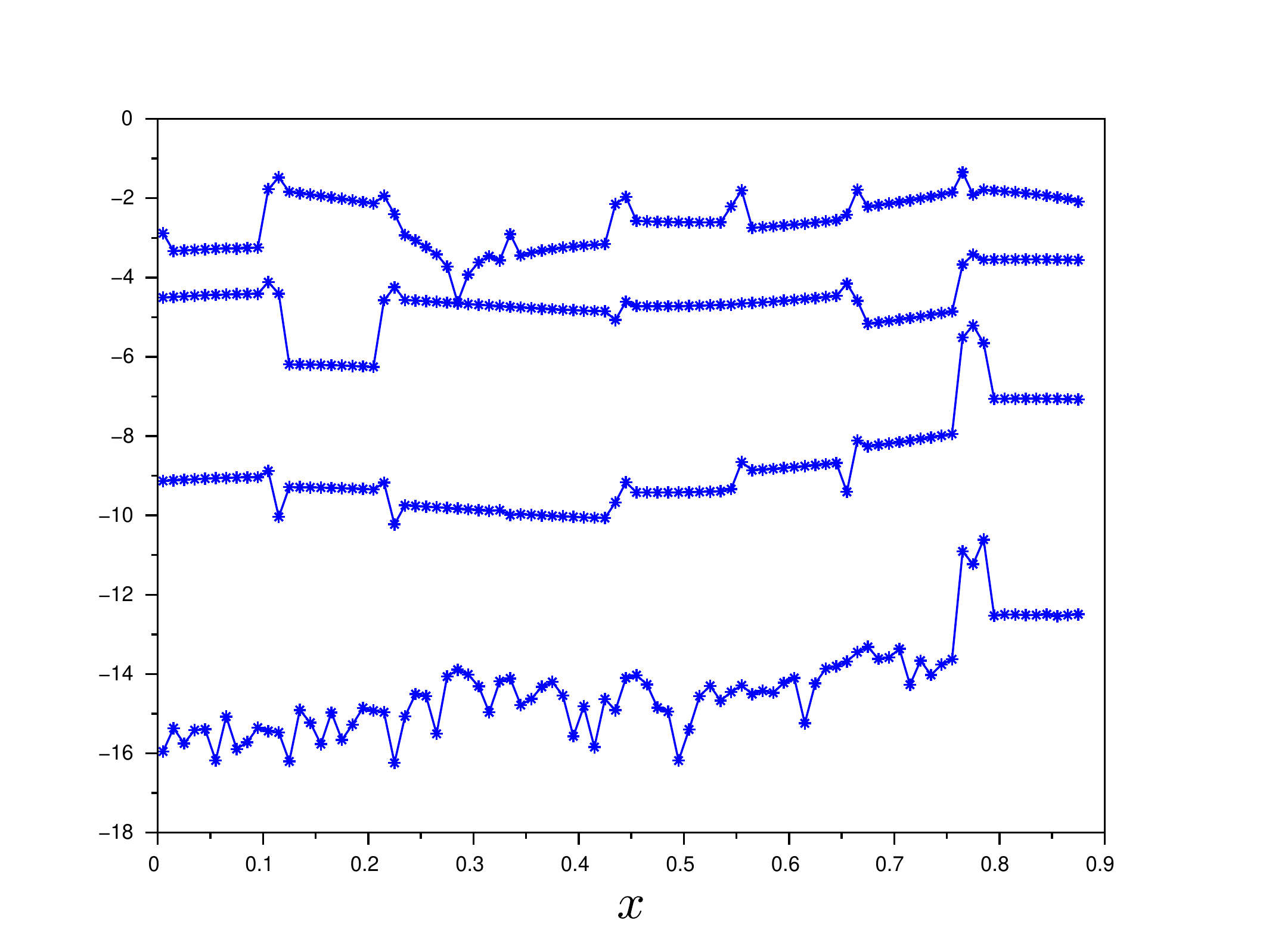}
  \caption{Illustration of the residual when RAS is used as a
    nonlinear solver (left), or as a preconditioner for
    Newton's method (right).}
 \label{Fig0}
\end{figure}

Like in the case of linear problems, where one usually accelerates the
Schwarz method, which is a fixed point iteration, using a Krylov
method, we can accelerate the non-linear fixed point iteration
(\ref{ClassicalSchwarz1dmodel}) using Newton's method. To do so, we
introduce two solution operators for the non-linear subdomain problems
in (\ref{ClassicalSchwarz1dmodel}),
\begin{equation}\label{SolOpSimple}
  u_1^n=G_1(u^{n-1}),\qquad  u_2^n=G_2(u^{n-1}),
\end{equation}
with which the classical parallel Schwarz method
(\ref{ClassicalSchwarz1dmodel}) can now be written in compact form,
even for many subdomains $i=1,\cdots,I$, as
\begin{equation}\label{SchwarzIterCompact}
  u^{n}=\sum_{i=1}^I{\widetilde P}_iG_i(u^{n-1})=:{\cal G}_1(u^{n-1}).
\end{equation}
As shown in the introduction, this fixed point iteration can be used
as a preconditioner for Newton's method, which means to apply Newton's
method to the non-linear equation
\begin{equation}\label{OneLevel}
  \tilde{\cal F}_1(u):={\cal G}_1(u)-u=\sum_{i=1}^I{\widetilde P}_iG_i(u)-u=0,
\end{equation}
because it is this equation that holds at the fixed point of iteration
(\ref{SchwarzIterCompact}). We call this method one level RASPEN
(Restricted Additive Schwarz Preconditioned Exact Newton). We show in
Figure \ref{Fig0} as an example the residual of the nonlinear RAS
iterations and using RASPEN as a preconditioner for Newton
when solving the Forchheimer equation with 8 subdomains from the
numerical section.
We observe that the residual of the non-linear RAS method is
concentrated at the interfaces, since it must be zero inside the
subdomains by construction.  Thus, when Newton's method is used to
solve \eqref{OneLevel}, it only needs to concentrate on reducing the
residual on a small number of interface variables.  This explains the
fast convergence of RASPEN shown on the right of Figure \ref{Fig0},
despite the slow convergence of the underlying RAS iteration.

Suppose we also want to include a coarse grid correction step in the
Schwarz iteration (\ref{ClassicalSchwarz1dmodel}), or equivalently in
(\ref{SchwarzIterCompact}). Since the problem is non-linear, we need
to use the Full Approximation Scheme (FAS) from multigrid to do so,
see for example \cite{briggs2000multigrid,henson2003multigrid}: given
an approximate solution $u^{n-1}$, we compute the correction $c$ by
solving the non-linear coarse problem
\begin{equation}\label{FAS}
  {\cal L}^c(R_0u^{n-1}+c)={\cal L}^c(R_0u^{n-1})+\tilde{R}_0(f-{\cal L}(u^{n-1})),
\end{equation}
where ${\cal L}^c$ is a coarse approximation of the non-linear problem
(\ref{1dmodel}) and $R_0$ is a restriction operator. This correction
$c:=C_0(u^{n-1})$ is then added to the iterate to get the new corrected value
\begin{equation}\label{coarsecorradd}
  u^{n-1}_{new}=u^{n-1}+P_0C_0(u^{n-1}),
\end{equation}
where $P_0$ is a suitable prolongation operator. Introducing this new
approximation from (\ref{coarsecorradd}) at step $n-1$ into the
subdomain iteration formula (\ref{SchwarzIterCompact}), we obtain the
method with integrated coarse correction
\begin{equation}
  u^{n}=\sum_{i=1}^I{\widetilde
    P}_iG_i(u^{n-1}+P_0C_0(u^{n-1}))=:{\cal G}_2(u^{n-1}).
\end{equation}
This stationary fixed point iteration can also be accelerated using
Newton's method: we can use Newton to solve the non-linear equation
\begin{equation}\label{TwoLevel}
  \tilde{\cal F}_2(u):={\cal G}_2(u)-u=\sum_{i=1}^I{\widetilde P}_iG_i(u+P_0C_0(u))-u=0.
\end{equation}
We call this method two level FAS-RASPEN.\\

We have written the coarse step as a correction, but not the
subdomain steps. This can however also be done, by simply rewriting
(\ref{SchwarzIterCompact}) to add and subtract the previous iterate,
\begin{equation}\label{SchwarzIterCompactC}
  u^{n}=\smash[t]{u^{n-1}+\sum_{i=1}^I{\widetilde P}_i\underbrace{(G_i(u^{n-1})-R_iu^{n-1})}_{=:C_i(u^{n-1})}
    =u^{n-1}+\sum_{i=1}^I{\widetilde P}_iC_i(u^{n-1}),}
\end{equation}
where we have assumed that $\sum_i \tilde{P}_iR_i = I_V$, the identity on the vector space, see Assumption 1 in the next section. Together with the coarse grid correction (\ref{coarsecorradd}), this
iteration then becomes
\begin{equation}\label{SchwarzIterCompactCC}
  u^{n}=\smash[t]{u^{n-1}+P_0C_0(u^{n-1})+\sum_{i=1}^I{\widetilde P}_iC_i(u^{n-1}+P_0C_0(u^{n-1})),}
\end{equation}
which can be accelerated by solving with Newton the equation
\begin{equation}\label{TwoLevelC}
  \tilde{\cal F}_2(u):=P_0C_0(u)+\sum_{i=1}^I{\widetilde P}_iC_i(u+P_0C_0(u))=0.
\end{equation}
This is equivalent to $\tilde{\cal F}_2(u)=0$ from (\ref{TwoLevel}), only
written in correction form.

\section{Definition of RASPEN and Comparison with ASPIN}\label{ComparisonSec}

We now define formally the one- and two-level versions of the RASPEN method and compare it with
the respective ASPIN methods. We consider a non-linear function $F:
V\rightarrow V'$, where $V$ is a Hilbert space, and the
non-linear problem of finding $u\in V$ such that
\begin{equation}\label{eqtobesolved}
  F(u) = 0.
\end{equation}
Let $V_i,i=1,\ldots,I$ be Hilbert spaces, which would generally be
subspaces of $V$. We consider for all $i=1,\ldots,I$ the linear
restriction and prolongation operators $R_i : V\rightarrow V_i$,
$P_i: V_i \rightarrow V$, 
as well as the ``restricted'' prolongation ${\widetilde P}_i: V_i \rightarrow V$.\\

\begin{assumption}
We assume that $R_i$  and $P_i$ satisfy for $i=1,\ldots,I$
$$
  R_iP_i = I_{V_i},\quad \mbox{the identity on $V_i$,}
$$
and that $R_i$  and ${\widetilde P}_i$ satisfy
$
\displaystyle  \;\;\sum_{i=1}^I {\widetilde P}_i R_i = I_V.
$\\
\end{assumption}

These are all the assumptions we need in what follows, but it is
helpful to think of the restriction operators $R_i$ as classical
selection matrices which pick unknowns corresponding to the subdomains
$\Omega_i$, of the prolongations $P_i$ as $R_i^T$, and of the
${\widetilde P}_i$ as extensions based on a non-overlapping decomposition.

\subsection{One- and two-level RASPEN}

We can now formulate precisely the RASPEN method from the previous
section: we define the local inverse $G_i: V\rightarrow V_i$ to be
solutions of
\begin{equation}\label{GDef}
  R_i F\left(P_i G_i(u) + (I-P_iR_i)u\right) = 0.
\end{equation}
In the usual PDE framework, this corresponds to solving locally on the
subdomain $i$ the PDE problem on $V_i$ with Dirichlet boundary condition
given by $u$ outside of the subdomain $i$, see (\ref{SolOpSimple}).
Then, one level RASPEN solves the non-linear equation
\begin{equation}\label{RASPEN}
  \tilde{\cal F}_1(u)  = \sum_{i=1}^I {\widetilde P}_i G_i(u)-u = 0,
\end{equation}
using Newton's method, see (\ref{OneLevel}). The preconditioned nonlinear
function (\ref{RASPEN}) corresponds to the fixed point iteration
\begin{equation}\label{RAS}
  u^{n}=\smash[t]{\sum_{i=1}^I{\widetilde P}_iG_i(u^{n-1}),}
\end{equation}
see (\ref{SchwarzIterCompact}). Equivalently, the RASPEN equation \eqref{RASPEN}
can be written in correction form as
\begin{equation}\label{RASPEN_corr}
\tilde{\cal F}_1(u)  = \smash[t]{\sum_{i=1}^I {\widetilde P}_i (G_i(u)-R_iu) =: \sum_{i=1}^I {\widetilde P}_iC_i(u),}
\end{equation}
where we define the corrections $C_i(u):= G_i(u) - R_iu$. This way, the subdomain solves
\eqref{GDef} can be written in terms of $C_i(u)$ as
\begin{equation}\label{Cdef}
  R_i F(u+P_iC_i(u)) = 0.
\end{equation}
In the special case where $F(u) = Au - b$ is affine, \eqref{Cdef} reduces to
$$ R_iA(u + P_iC_i(u)) - R_ib = 0 \implies C_i(u) = A_i^{-1}R_i(b-Au), $$
where $A_i = R_iAP_i$ is the subdomain matrix. This implies
$$ \tilde{\cal F}_1(u) = \sum_{i=1}^I {\widetilde P}_iA_i^{-1}R_i(b-Au), $$
and we immediately see that the Jacobian is the matrix $A$ preconditioned by the restricted additive
Schwarz (RAS) preconditioner $\sum_{i=1}^I {\widetilde P}_iA_i^{-1}R_i$. Thus, if a Krylov method is used to solve the outer system, our method is equivalent to the Krylov-accelerated one-level RAS method in the linear case.

To define the two-level variant, we introduce a coarse space $V_0$ and the linear restriction and prolongation operators $R_0 : V\rightarrow V_0$, $P_0 : V_0\rightarrow V$.
Let $F_0: V_0 \rightarrow V'_0$ be the coarse non-linear function,
which could be defined by using a coarse discretization of the
underlying problem, or using a Galerkin approach we use here, namely
\begin{equation}\label{CoarseProblem}
  F_0(u_0) =  {\widetilde R}_0 F(P_0(u)).
\end{equation}
Here, ${\widetilde R}_0:V' \rightarrow V'_0$ is a projection operator
that plays the same role as $R_0$, but in the residual space.
In two-level FAS-RASPEN, we use the well established
non-linear coarse correction $C_0(u)$ from the full approximation
scheme already shown in (\ref{FAS}), which in the rigorous context
of this section is defined by
\begin{equation}\label{FASCC}
  F_0(C_0(u) + R_0 u) = F_0(R_0 u) -  {\widetilde R}_0 F(u).
\end{equation}
This coarse correction is used in a multiplicative fashion in RASPEN,
i.e. we solve with Newton the preconditioned non-linear system
\begin{equation}\label{FASRASPEN}
  \widetilde{\cal F}_2(u) = P_0 C_0(u)+ \sum_{i=1}^n {\widetilde P}_i C_i(u + P_0C_0(u)) = 0.
\end{equation}
This corresponds to the non-linear two-level fixed point iteration
$$
  u^{n+1} =u^n+P_0 C_0(u^n)+\sum_{i=1}^n {\widetilde P}_i C_i(u^n + P_0C_0(u^n)),
$$
with $C_0(u^n)$ defined in (\ref{FASCC}) and $C_i(u^n)$ defined in
(\ref{Cdef}). This iteration is convergent, as we can see in the next section in Figure
\ref{Fig1} in the right column.
In the special case of an affine residual function $F(u) = Au - b$, a simple calculation shows that
$$\widetilde{\cal F}_2(u) =
\left(P_0A_0^{-1}{\widetilde R}_0 + \sum_{i=1}^I {\widetilde P}_iA_i^{-1}R_i(I_V - P_0A_0^{-1}{\widetilde R}_0)\right)(b - Au), $$
where we assumed that the coarse function $F_0 = A_0u_0 - b_0$ is also linear. Thus, \FK{ in the linear case}, two-level RASPEN
corresponds to preconditioning by a two-level RAS preconditioner, where the coarse grid correction is applied multiplicatively.

\subsection{Comparison of one-level variants}\label{OneLevelSection}
In order to compare RASPEN with the existing ASPIN method, we recall the precise definition of one-level ASPIN from
\cite{cai2002nonlinearly}, which gives a different reformulation
${\cal F}_1(u) = 0$ of the original equation (\ref{eqtobesolved}) to be
solved.  In ASPIN, one also defines for $u\in V$ and for all $i=1,\cdots,I$
the corrections as in \eqref{Cdef}, i.e., we define $C_i(u)\in V_i$ such that
\begin{equation*}
  R_i F(u+P_iC_i(u)) = 0,
\end{equation*}
where $P_iC_i(u)$ are called $T_i$ in \cite{cai2002nonlinearly}.
Then, the one-level ASPIN preconditioned function is defined by
\begin{equation}\label{ASPIN}
  {\cal F}_1(u) = \sum_{i=1}^I P_i C_i(u),
\end{equation}
and the preconditioned system ${\cal F}_1(u) = 0$ is solved using a
Newton algorithm with an inexact Jacobian, see Section \ref{JacobianSection}. The ASPIN preconditioner
also has a corresponding fixed point iteration: adding and subtracting
$P_i R_i u$ in the definition (\ref{Cdef}) of the corrections $C_i$,
we obtain
$$
  R_i F(u+P_iC_i(u)) = R_i F\left(P_i (R_i u+C_i(u)) + u-P_i R_i u\right)=0,
$$
which implies, by comparing with (\ref{GDef}) and assuming existence and
uniqueness of the solution to the local problems, that
$$
  G_i(u) = R_iu+C_i(u).
$$
We therefore obtain for one-level ASPIN
\begin{equation}\label{ASPINS}
  {\cal F}_1(u)=\sum_{i=1}^I P_i C_i(u)
  = \sum_{i=1}^I P_i G_i(u)-\sum_{i=1}^I P_iR_i u,
\end{equation}
which corresponds to the non-linear fixed point iteration
\begin{equation}\label{AS}
  u^{n}=u^{n-1}+\sum_{i=1}^I P_i C_i(u^{n-1})=u^{n-1}
    -\sum_{i=1}^I P_iR_i u^{n-1} + \sum_{i=1}^I P_i G_i(u^{n-1}).
\end{equation}
This iteration is not convergent in the overlap, already in the linear
case, see \cite{Efstathiou:2003:WRA,gander2008schwarz}, and needs a
relaxation parameter to yield convergence, see for example
\cite{dryja1997nonlinear} for the non-linear case. This can be seen
directly from (\ref{AS}): if an overlapping region belongs to $K$
subdomains, then the current iterate $u^n$ is subtracted $K$ times
there, and then the sum of the $K$ respective subdomain solutions are
added to the result. This redundancy is avoided in our formulation
\eqref{RAS}. The only interest in using an additive correction in the
overlap is that in the linear case, the preconditioner remains
symmetric for a symmetric problem.

\begin{figure}
  \centering
  \mbox{\includegraphics[width=0.49\textwidth]{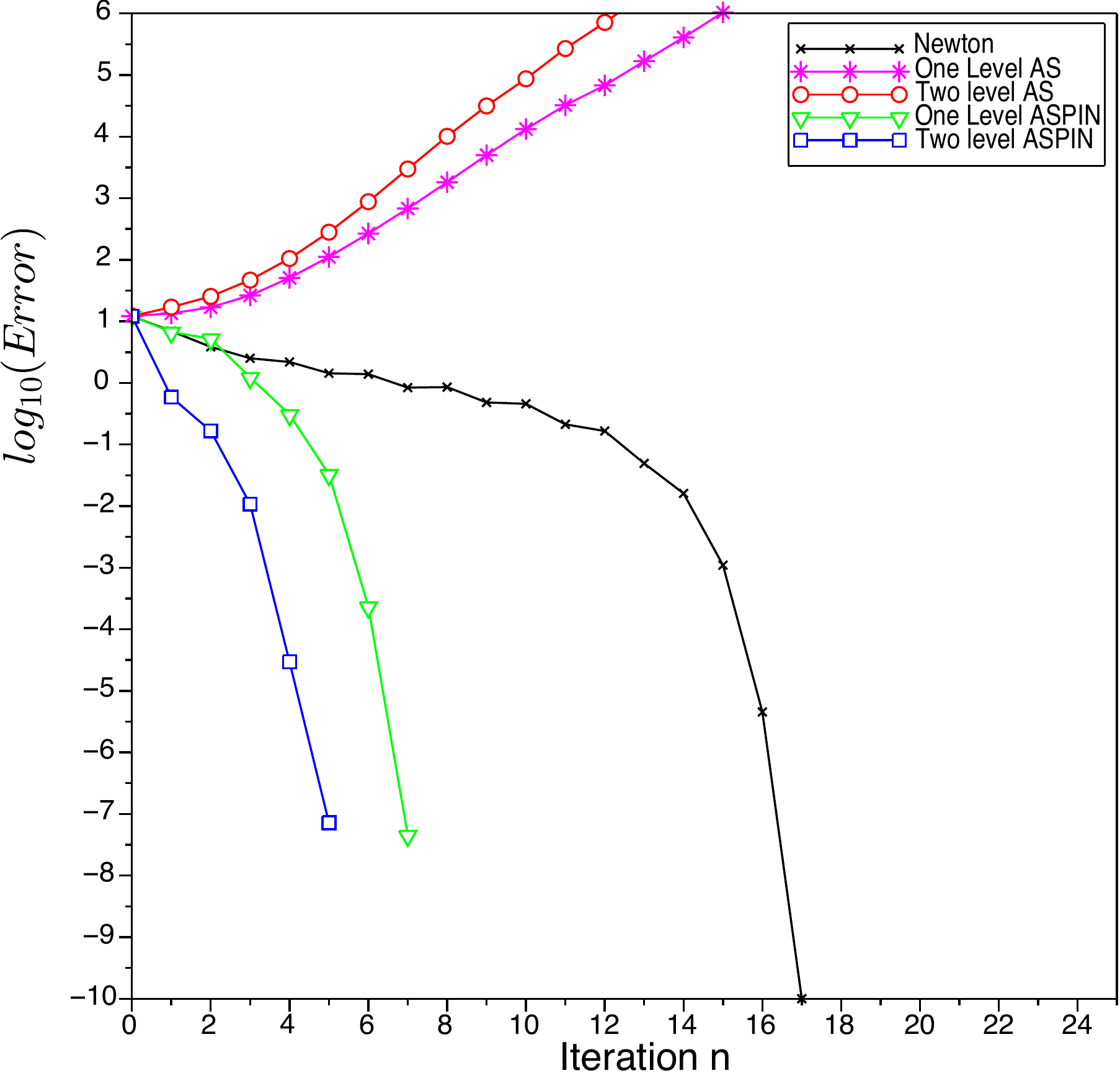}
  \includegraphics[width=0.49\textwidth]{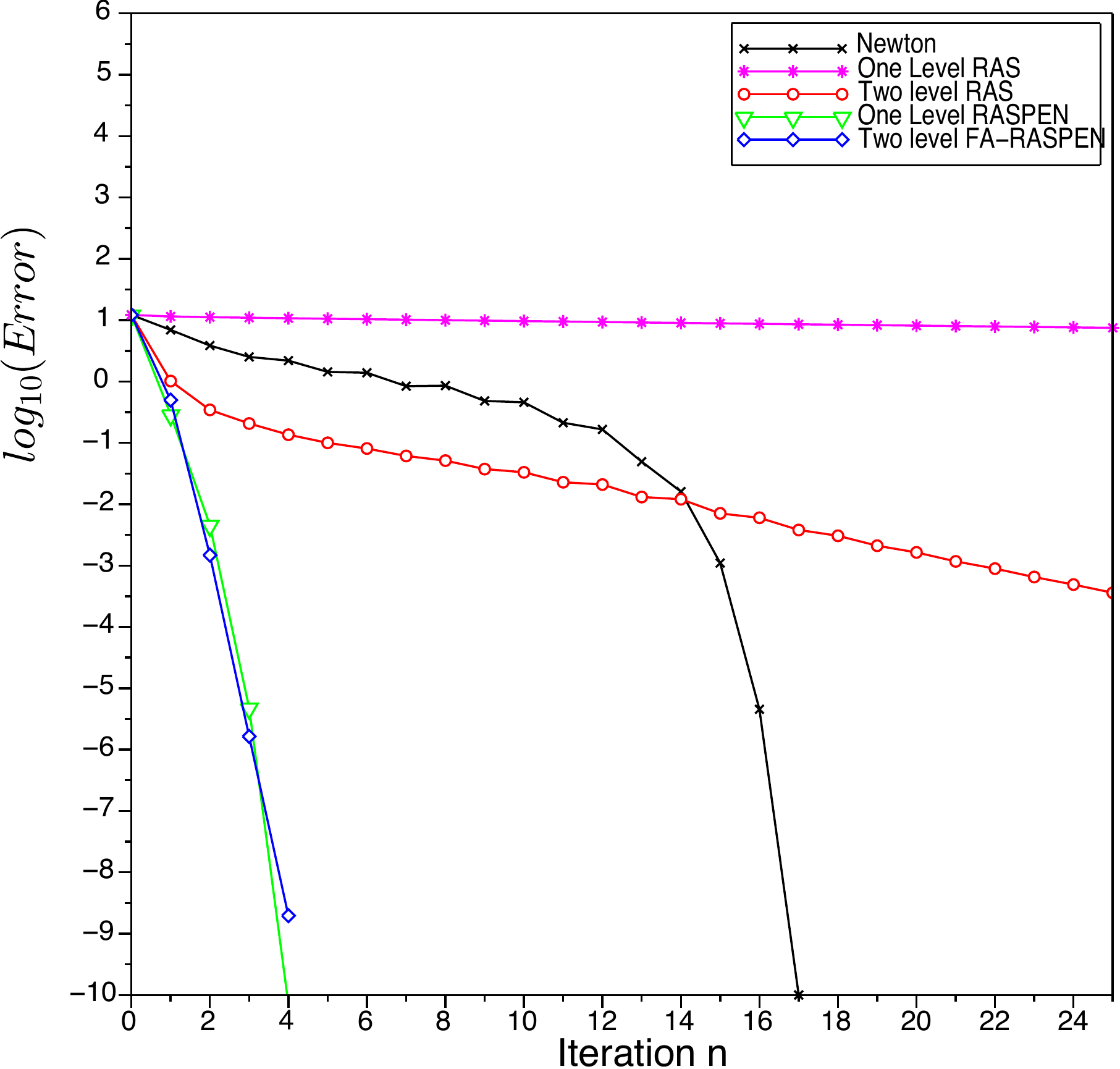}}
  \mbox{\includegraphics[width=0.49\textwidth]{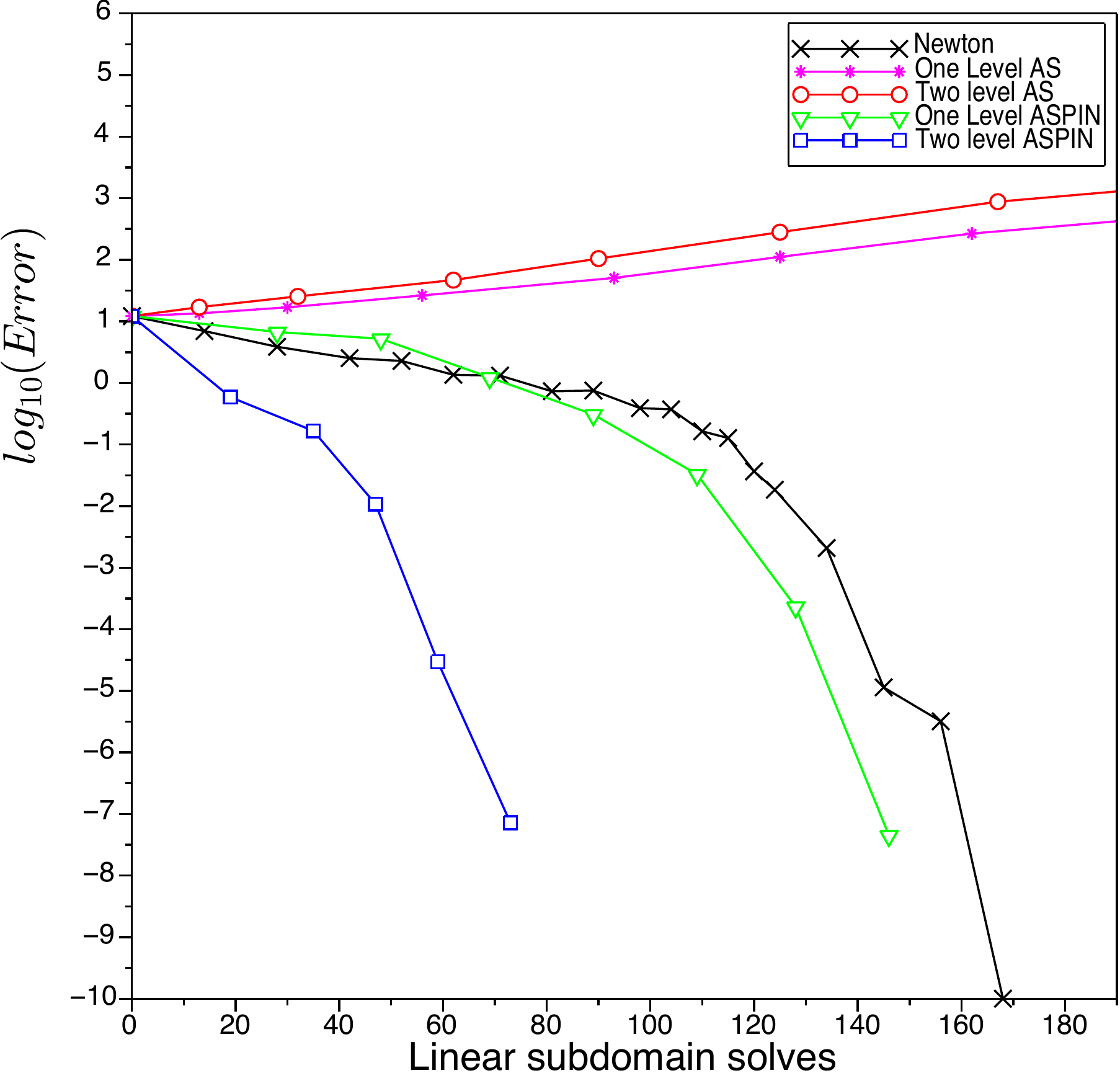}
  \includegraphics[width=0.49\textwidth]{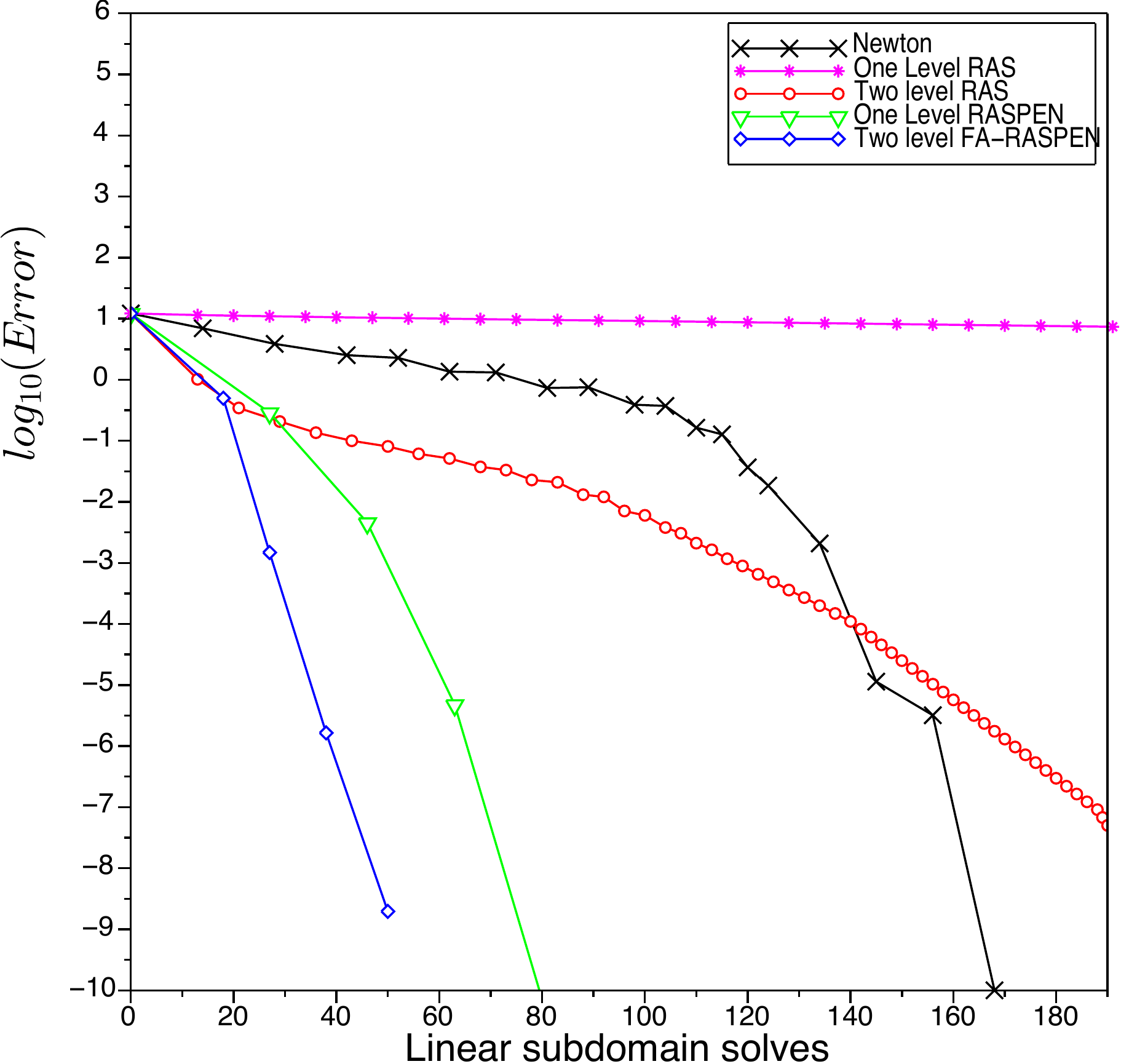}}
  \caption{Error as function of non-linear iteration numbers in the top
    row, and as number of subdomain solves in the bottom row,
    for ASPIN (left), and RASPEN (right).}
  \label{Fig1}
\end{figure}

We show in Figure \ref{Fig1} a numerical comparison of the two
methods, together with Newton's method applied directly to the
non-linear problem, for the first example of the Forchheimer
equation from Section \ref{ForchSec} on a domain of unit size with 8
subdomains, overlap $3h$, with $h=1/100$. In these comparisons, we use
ASPIN first as a fixed-point iterative solver (labelled AS for
Additive Schwarz), and then as a preconditioner. We do the same for
our new nonlinear iterative method, which in the figures are labelled
RAS for Restricted Additive Schwarz.
We see from this numerical experiment that ASPIN as an iterative
solver (AS) does not converge, whereas RASPEN used as an iterative
solver (RAS) does, both with and without coarse grid. Also note that
two-level RAS is faster than Newton directly applied to the non-linear problem for
small iteration counts, before the superlinear convergence of Newton kicks in.
The fact that RASPEN is based on a convergent iteration, but not ASPIN, has an
important influence also on the Newton iterations when the methods are used as preconditioners: the ASPIN preconditioner requires more Newton
iterations to converge than RASPEN does. At first sight, it might be surprising that in
RASPEN, the number of Newton iterations with and without coarse grid is
almost the same, while ASPIN needs more iterations without coarse
grid. In contrast to the linear case with Krylov acceleration, it is
not the number of Newton iterations that depends on the number
of subdomains, but the number of linear inner iterations within
Newton, which grows when no coarse grid is present. We show this in
the second row of Figure \ref{Fig1}, where now the error is plotted as a
function of the maximum number of linear subdomain solves used in each
Newton step, see Subsection \ref{subseconelevelnum}. With this more
realistic measure of work, we see that both RASPEN and ASPIN converge
substantially better with a coarse grid, but RASPEN needs much fewer subdomain solves than ASPIN does.

\subsection{Comparison of two-level variants}

We now compare two-level FAS-RASPEN with the two-level ASPIN method
of \cite{marcinkowski2005parallel}. Recall that the two-level FAS-RASPEN
consists of applying Newton's method to \eqref{FASRASPEN},
\begin{equation*}
  \widetilde{\cal F}_2(u) = P_0 C_0(u)+ \sum_{i=1}^n {\widetilde P}_i C_i(u + P_0C_0(u)) = 0,
\end{equation*}
where the corrections $C_0(u)$ and $C_i(u)$ are defined in \eqref{FASCC} and
\eqref{Cdef} respectively. Unlike FAS-RASPEN, two-level ASPIN requires
the solution $u^*_0\in V_0$ to the coarse problem, i.e., $F_0(u^*_0) = 0$,
which can be computed in a preprocessing step.

In two-level ASPIN, the coarse correction $C_0^A: V\rightarrow V_0$ is
defined by
\begin{equation}\label{ASPINCC}
  F_0(C_0^A(u) + u_0^*) = - {\widetilde R}_0 F(u),
\end{equation}
and the associated two-level ASPIN function uses the coarse correction
in an additive fashion, i.e. Newton's method is used to solve
\begin{equation}\label{2LevelASPIN}
  {\cal F}_2(u) = P_0 C_0^A(u) + \sum_{i=1}^I P_i C_i(u) = 0,
\end{equation}
with $C_0^A(u^n)$ defined in (\ref{ASPINCC}) and $C_i(u^n)$ defined in
(\ref{Cdef}).
This is in contrast to two-level FAS-RASPEN, where
the coarse correction $C_0(u)$ is computed from the well established
full approximation scheme, and is applied multiplicatively in
\eqref{FASRASPEN}. The fixed point iteration corresponding to
 \eqref{2LevelASPIN} is
$$
  u^{n+1} =u^n+P_0 C_0^A(u^n)+ \sum_{i=1}^I P_i C_i(u^n).
$$
Just like its one-level counterpart, two-level ASPIN is not convergent as a fixed-point iteration without a relaxation parameter,
see Figure \ref{Fig1} in the left column. Moreover, because the coarse correction is applied additively, the overlap between the coarse space and subdomains leads to slower convergence in the
Newton solver, which does not happen with FAS-RASPEN.

\subsection{Computation of Jacobian matrices}\label{JacobianSection}

When solving \eqref{RASPEN_corr}, \eqref{FASRASPEN}, \eqref{ASPINS}
and \eqref{2LevelASPIN} using Newton's method, one needs to repeatedly
solve linear systems involving Jacobians of the above functions. If
one uses a Krylov method such as GMRES to solve these linear systems,
like we do in this paper, then it suffices to have a procedure for
multiplying the Jacobian with an arbitrary vector.  In this section,
we derive the Jacobian matrices for both one-level and two-level
RASPEN in detail. We compare these expressions with ASPIN, which
approximates the exact Jacobian with an inexact one in an attempt to
reduce the computational cost, even though this can potentially slow
down the convergence of Newton's method. Finally, we show that this
approximation is not necessary in RASPEN: in fact, all the components
involved in building the exact Jacobian have already been computed
elsewhere in the algorithm, so there is little additional cost in
using the exact Jacobian compared with the approximate one.

\subsubsection{Computation of the one-level Jacobian matrices}
We now show how to compute the Jacobian matrices of ASPIN and RASPEN.
To simplify notation, we define
\begin{equation}
  u^{(i)}:= P_i G_i(u) + (I-P_iR_i)u\quad \mbox{and}\quad J(v) :=
    {dF\over du}(v)
\end{equation}
By differentiating (\ref{GDef}),
we obtain
\begin{equation}\label{dGdu}
{dG_i\over du}(u) = - (R_i J(u^{(i)}) P_i)^{-1}R_i J(u^{(i)}) + R_i.
\end{equation}
We deduce for the Jacobian of RASPEN from (\ref{RASPEN})
\begin{equation}\label{RASPENJac}
  {d\tilde{\cal F}_1 \over du}(u)
  =\sum_{i=1}^I {\widetilde P}_i {dG_i \over du}(u)-I
  =-\sum_{i=1}^I {\widetilde P}_i (R_i J(u^{(i)})
  P_i)^{-1}R_i J(u^{(i)}),
\end{equation}
since the identity cancels.  Similarly, we obtain for the Jacobian of
ASPEN (Additive Schwarz Preconditioned Exact Newton) in (\ref{ASPINS})
\begin{equation}\label{ASPENJac}
  {d{\cal F}_1 \over du}(u) =
    \sum_{i=1}^I P_i {dG_i \over du}(u)-\sum_{i=1}^I P_iR_i=
  -\sum_{i=1}^I P_i (R_i J(u^{(i)})
  P_i)^{-1}R_i J(u^{(i)}),
\end{equation}
since now the terms $\sum_{i=1}^I P_iR_i$ cancel.  In ASPIN, this
exact Jacobian is replaced by the inexact Jacobian
$$
  {d{\cal F}_1 \over du}^{inexact}(u)
    = -\left(\sum_{i=1}^I P_i (R_i J(u)P_i)^{-1}R_i\right) J(u).
$$
We see that this is equivalent to preconditioning the Jacobian $J(u)$
of $F(u)$ by the additive Schwarz
preconditioner, up to the minus sign. This can be convenient if one
has already a code for this, as it was noted in \cite{cai2002nonlinearly}.
The exact Jacobian is however also easily accessible,
since the Newton solver for the non-linear subdomain system
$R_i F\left(P_i G_i(u) + (I-P_iR_i)u\right) = 0$ already computes and factorizes the
local Jacobian matrix $R_i J(u^{(i)}) P_i$. Therefore, the only missing ingredient for computing the
exact Jacobian of ${\cal F}_1$ is the matrices $R_i J(u^{(i)})$, which only differ from $R_i J(u^{(i)}) P_i$ by a few additional columns, corresponding in the usual
PDEs framework to the derivative with respect to the Dirichlet boundary conditions.
In contrast, the computation of the inexact ASPIN Jacobian requires one to recompute the entire Jacobian
of $F(u)$ after the subdomain non-linear solves.

\subsubsection{Computation of the two-level Jacobian matrices}
We now compare the Jacobians for the two-level variants. For RASPEN,
we need to differentiate $\tilde{\cal F}_2$ with respect to $u$, where $\tilde{\cal F}_2$ is defined in
(\ref{FASRASPEN}):
\begin{equation*}
  \widetilde{\cal F}_2(u) = P_0 C_0(u)+ \sum_{i=1}^n {\widetilde P}_i C_i(u + P_0C_0(u)).
\end{equation*}
To do so, we need
$\frac{dC_0}{du}$ and $\frac{dC_i}{du}$ for $i=1,\ldots, I$. The former can be obtained by differentiating
\eqref{FASCC}:
$$ F_0'(R_0u + C_0(u)) \left(R_0 + \frac{dC_0}{du}\right) = F_0'(R_0u)R_0 -  {\widetilde R}_0 F'(u). $$
Thus, we have
\begin{equation}\label{eq:C0}
\frac{dC_0}{du} = -R_0 + \hat{J}_0^{-1}(J_0 R_0 -  {\widetilde R}_0 J(u)),
\end{equation}
where $J_0 = F_0'(R_0u)$ and $\hat{J}_0 = F_0'(R_0u + C_0(u))$.
Note that the two Jacobian matrices are evaluated at different arguments, so no cancellation is possible in \eqref{eq:C0} except in special cases (e.g., if $F_0$ is an affine function). Nonetheless, they are readily available: $\hat{J}_0$ is simply the Jacobian for the non-linear coarse solve, so it would have already been calculated and factorized by Newton's method. $J_0$ would also have been calculated during the coarse Newton iteration if $R_0u$ is used as the initial guess.

We also need $\frac{dC_i}{du}$
from the subdomain solves. From the relation
$G_i(u) = R_iu + C_i(u), $
we deduce immediately from \eqref{dGdu} that
\begin{equation}\label{eq:Ci}
\frac{dC_i}{du} = \frac{dG_i}{du} - R_i = -(R_iJ(u^{(i)})P_i)^{-1}R_iJ(u^{(i)}),
\end{equation}
where $u^{(i)} = u + P_iC_i(u). $
Thus, the Jacobian for the two-level RASPEN function is
\begin{equation}\label{dF2tdu}
\frac{d\widetilde{\mathcal{F}}_2}{du} = P_0 \frac{dC_0}{du} - \sum_i \tilde P_i (R_iJ(v^{(i)})P_i)^{-1}R_i
J(v^{(i)})\left(I + P_0\frac{dC_0}{du}\right),
\end{equation}
where $\frac{dC_0}{du}$ is given by \eqref{eq:C0} and
$v^{(i)} = u + P_0C_0(u) + P_iC_i(u + P_0C_0(u)). $

For completeness, we compute the Jacobian for two-level ASPIN. First, we obtain $\frac{dC_0^A}{du}$ by differentiating \eqref{ASPINCC},
which gives
\begin{equation}\label{eq:C0_aspin}
\frac{dC_0^A}{du} = - \hat{J}_0^{-1} {\widetilde R}_0 J(u),
\end{equation}
where $\hat{J}_0 = F_0'(C_0^A(u)+u_0^*).$
In addition, two-level ASPIN uses as approximation for \eqref{eq:Ci}
\begin{equation}\label{eq:Ci_aspin}
\frac{dC_i}{du} \approx -(R_iJ(u)P_i)^{-1}R_i J(u).
\end{equation}
Thus, the inexact Jacobian for the two-level ASPIN function is
\begin{equation}\label{dF2du}
\frac{d\mathcal{F}_2}{du} \approx -P_0
\hat{J}_0^{-1} {\widetilde R}_0 J(u)  - \sum_i P_i (R_iJ(u)P_i)^{-1}R_i
J(u).
\end{equation}
Comparing \eqref{dF2tdu} with \eqref{dF2du}, we see two major differences. First, $dC_0/du$ only simplifies to $-{\widetilde R}_0J(u)$ if $J_0 = \hat{J}_0$, i.e., if $F_0$ is affine. Second, \eqref{dF2tdu} resembles a two-stage multiplicative preconditioner, whereas \eqref{dF2du} is of the additive type. This is due to the fact that the coarse correction in two-level RASPEN is applied multiplicatively, whereas two-level ASPIN uses an additive correction.

\section{Numerical experiments}\label{NumSec}

In this section, we compare the new non-linear preconditioner RASPEN to ASPIN for the
Forchheimer model, which generalizes the linear Darcy model in porous media flow \cite{forchheimer1901, ward1964, chen2006}, and for a 2D non-linear diffusion problem that appears in \cite{fenics}.

\subsection{Forchheimer model and discretization}\label{ForchSec}

Let us consider the Forchheimer parameter $\beta > 0$,
the permeability $\lambda\in L^\infty(\Omega)$ such that $0 < \lambda_{min} \leq \lambda(x) \leq
\lambda_{max}$ for all $x\in \Omega$, and the function
$q(g) = \mbox{sgn}(g) {- 1 + \sqrt{1+ 4\beta |g|}\over  2\beta}$.
The Forchheimer model on the interval $\Omega = (0,L)$ is defined by the equation
\begin{equation}\label{forch}
\left\{\begin{array}{rcll}
\left( q(-\lambda(x) u(x)')\right)' &=& f(x) &\mbox{ in } \Omega,\\
u(0) &=& u^D_0, \\
u(L) &=& u^D_L.
\end{array}\right.
\end{equation}
Note that at the limit when $\beta \rightarrow 0^+$, we recover the linear Darcy equation.
We consider a 1D mesh defined by the $M+1$ points
$$
0=x_{1\over 2} < \cdots  < x_{K + {1\over 2}} < \cdots < x_{M+{1\over 2}} = L.
$$
The cells are defined by $K = (x_{K-{1\over 2}}, x_{K+{1\over 2}})$ for
$K\in{\cal M}=\{1,\cdots,M\}$ and their center by $x_K = {x_{K-{1\over 2}} + x_{K+{1\over 2}} \over 2}$.
The Forchheimer model (\ref{forch}) is discretized using a
Two Point Flux Approximation (TPFA) finite volume scheme.
We define the TPFA transmissibilities by
$$
\left\{\begin{array}{r@{\,\,}c@{\,\,}l}
&&\displaystyle T_{K+{1\over 2}}  = {1
\over {|x_{K+{1\over 2}}-x_K| \over \lambda_K} + {|x_{K+1} - x_{K+{1\over 2}}| \over \lambda_{K+1} } } \mbox{ for } K=1,\cdots,M-1\\
&&\displaystyle T_{1\over 2} = {\lambda_{1} \over |x_1 -x_{1\over 2}|},\quad
\displaystyle T_{M+ {1\over 2}} =  {\lambda_{M} \over |x_{M + {1\over 2}}-x_M|},
\end{array}\right.
$$
with  $\lambda_K = {1\over |x_{K+{1\over 2}}- x_{K-{1\over 2}}| } \int_{x_{K-{1\over 2}}}^{x_{K+{1\over 2}}} \lambda(x)dx$.
Then, the $M$ cell unknowns $u_K$, $K\in {\cal M}$, are the solution of
the set of $M$ conservation equations
$$
\left\{\begin{array}{r@{\,\,}c@{\,\,}l}
\displaystyle q(T_{K+{1\over 2}}(u_{K}- u_{K+1})) +  q(T_{K-{1\over 2}}(u_{K}- u_{K-1})) &=& f_K,\,  K=2,\cdots,M-1\\
\displaystyle q(T_{{3\over 2}}(u_{1}- u_{2})) +  q(T_{{1\over 2}}(u_{1}- u_{0}^D)) &=& f_1,  \\
\displaystyle q(T_{M+{1\over 2}}(u_{M}- u_{L}^D)) +  q(T_{M-{1\over 2}}(u_{M}- u_{M-1})) &=& f_M,
\end{array}\right.
$$
with $f_K = \int_{x_{K-{1\over 2}}}^{x_{K+{1\over 2}}} f(x)dx$.
In the following numerical tests we will consider a uniform mesh
of cell size denoted by $h = {L\over M}$.

\subsubsection{One level variants}\label{subseconelevelnum}

We start from a non-overlapping decomposition of the set of cells
$$
\widetilde {\cal M}_i, i=1,\cdots,I,
$$
such that ${\cal M} = \bigcup_{i=1,\cdots,I} \widetilde {\cal M}_i$ and $\widetilde {\cal M}_i\cap \widetilde {\cal M}_j =\emptyset$
for all $i\neq j$.

The overlapping decomposition ${\cal M}_i, i=1,\cdots,I$ of the set of
cells is obtained by adding $k$ layers of cells to each $\widetilde
{\cal M}_i$ to generate overlap with the two neighbouring subdomains
$\widetilde {\cal M}_{i-1}$ (if $i>1$) and $\widetilde {\cal M}_{i+1}$
(if $i<I$) in the simple case of our one dimensional domain.

In the ASPIN framework, we set $V= {\Bbb R}^{\#{\cal M}}$, and $V_i=
{\Bbb R}^{\#{\cal M}_i}$, $i=1,\cdots,I$. The restriction
operators are then defined by
$$
(R_i v)_K = v_K \mbox{ for } K\in {\cal M}_i,
$$
and the prolongation operators are
$$
\left\{\begin{array}{r@{\,\,}c@{\,\,}l}
(P_i v_i)_K &=& v_K \mbox{ for } K\in {\cal M}_i,\\
(P_i v_i)_K &=& 0 \mbox{ for } K\not\in {\cal M}_i,
\end{array}\right.
\quad\text{ and }\quad
\left\{\begin{array}{r@{\,\,}c@{\,\,}l}
({\widetilde P}_i v_i)_K &=& v_K \mbox{ for } K\in \widetilde {\cal M}_i,\\
({\widetilde P}_i v_i)_K &=& 0 \mbox{ for } K\not\in \widetilde {\cal M}_i.
\end{array}\right.
$$
The coarse grid is obtained by the agglomeration of the cells in each
$\widetilde {\cal M}_i$ defining a coarse mesh of $(0,L)$.

Finally, we set $V_0 = {\Bbb R}^I$. In the finite volume framework, we
define for all $v\in V$
\begin{align*}
(R_0 v)_i &= {1\over \#\widetilde {M}_i} \sum_{K\in \widetilde {\cal M}_i} v_K \quad\mbox{for all } i=1,\cdots,I,\\
({\widetilde R}_0 v)_i &= \sum_{K\in \widetilde {\cal M}_i} v_K \;\;\quad \qquad \mbox{for all } i=1,\cdots,I.
\end{align*}
\FK{In our case of a uniform mesh, $R_0$ corresponds to the mean value in
the coarse cell $i$ for cellwise constant functions on ${\cal M}$, whereas
${\widetilde R}_0$ corresponds to the aggregate flux over the coarse cell ${\widetilde {\cal M}}_i$.}

For $v_0\in V_0$, its prolongation $v=P_0v_0\in V$ is obtained by the
piecewise linear interpolation $\varphi(x)$ on $(0,x_1)$,
$(x_1,x_2),\cdots,(x_I,L)$ where the $x_i$ are the centers of the
coarse cells, and $\varphi(x_i) = (v_0)_i$, $i=1,\cdots,I$,
$\varphi(0)=0$, $\varphi(L)=0$.  Then, $v=P_0v_0$ is defined by $v_K =
\varphi(x_K)$ for all $K\in {\cal M}$.  The coarse grid operator $F_0$
is defined by $F_0(v_0) = {\widetilde R}_0 F(P_0 v_0)$ for all $v_0\in
V_0$.

\begin{figure}
  \centering \mbox{
    \includegraphics[height=0.32\textwidth]{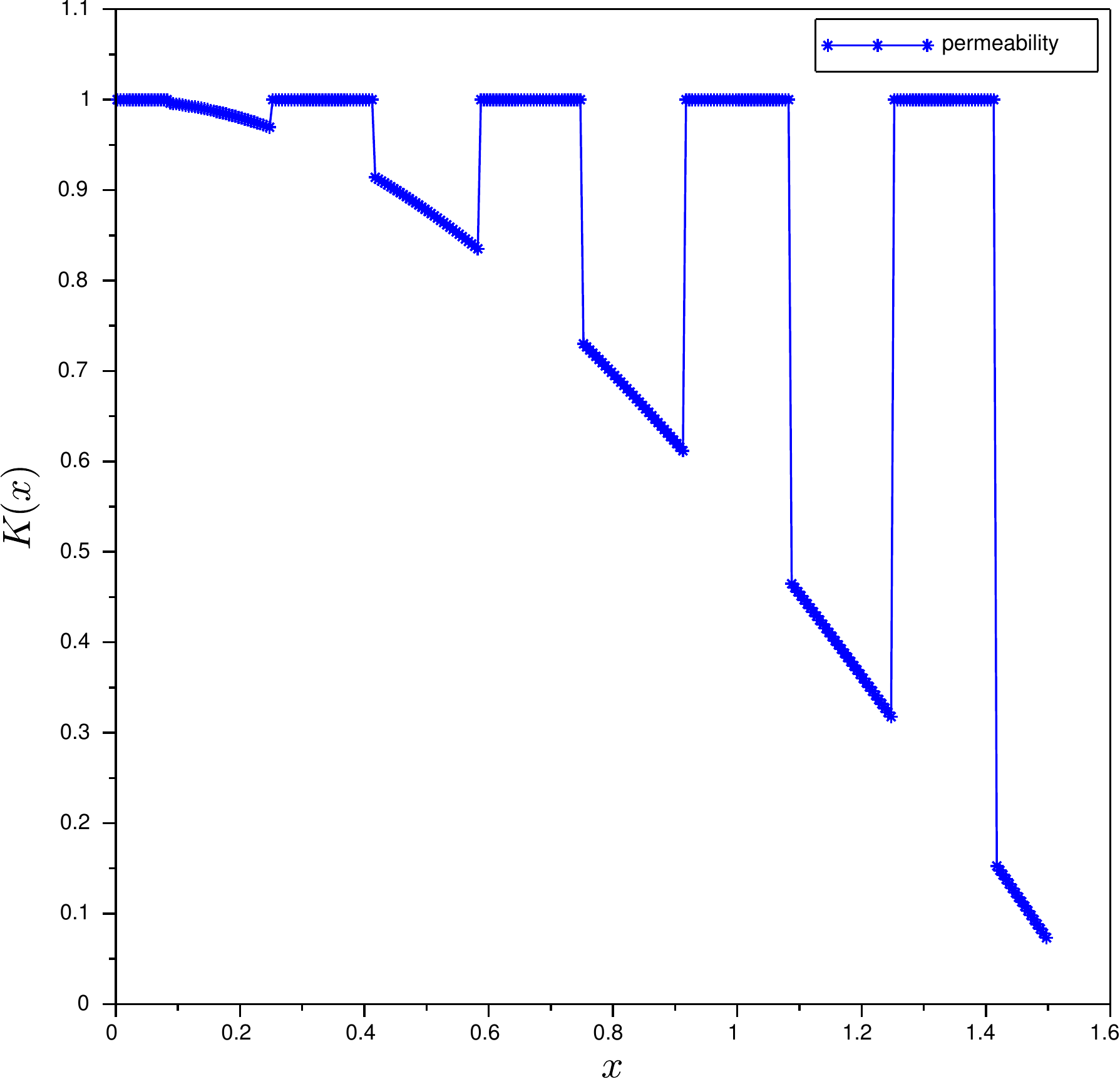}
    \includegraphics[height=0.32\textwidth]{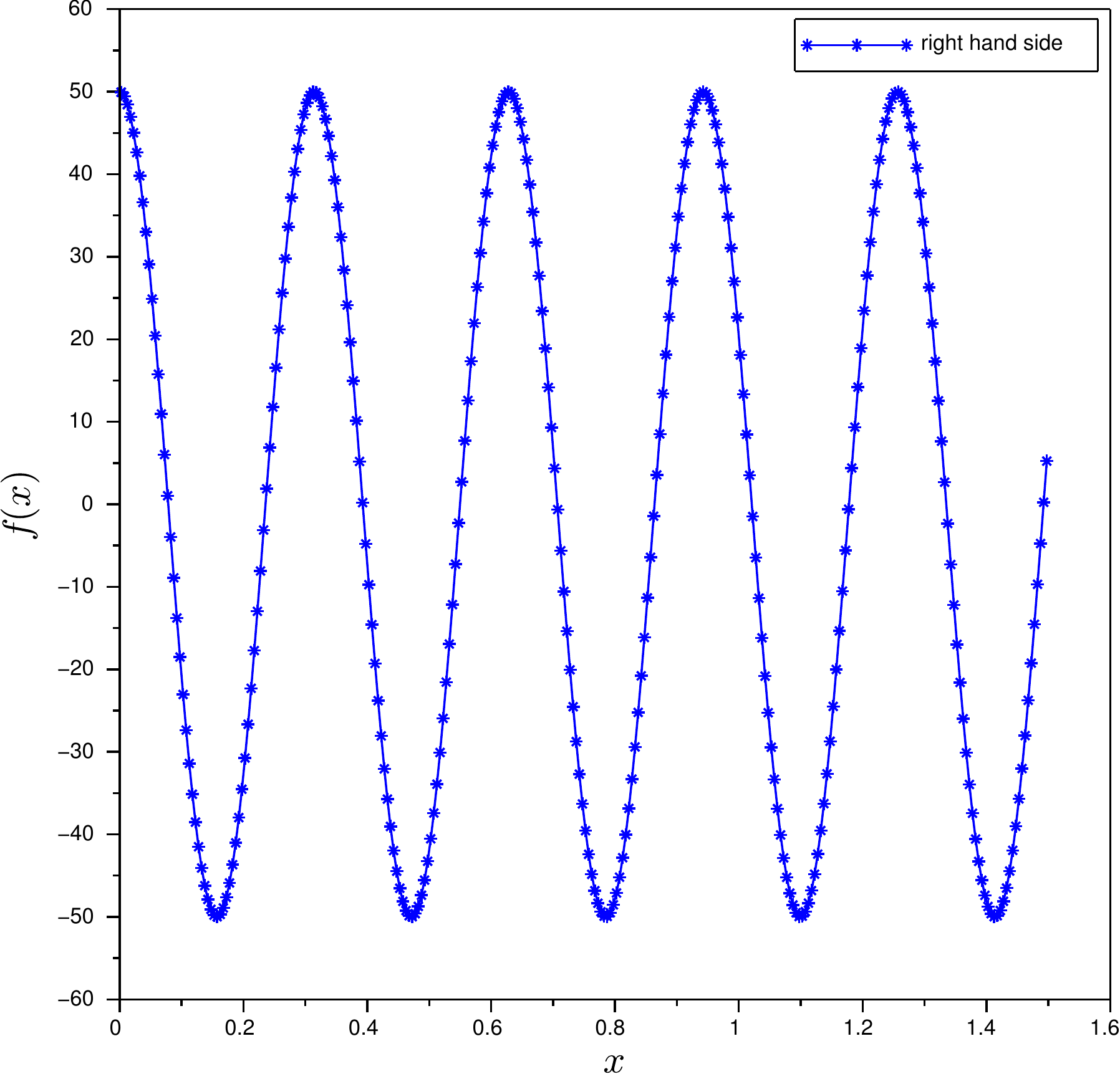}
    \includegraphics[height=0.32\textwidth]{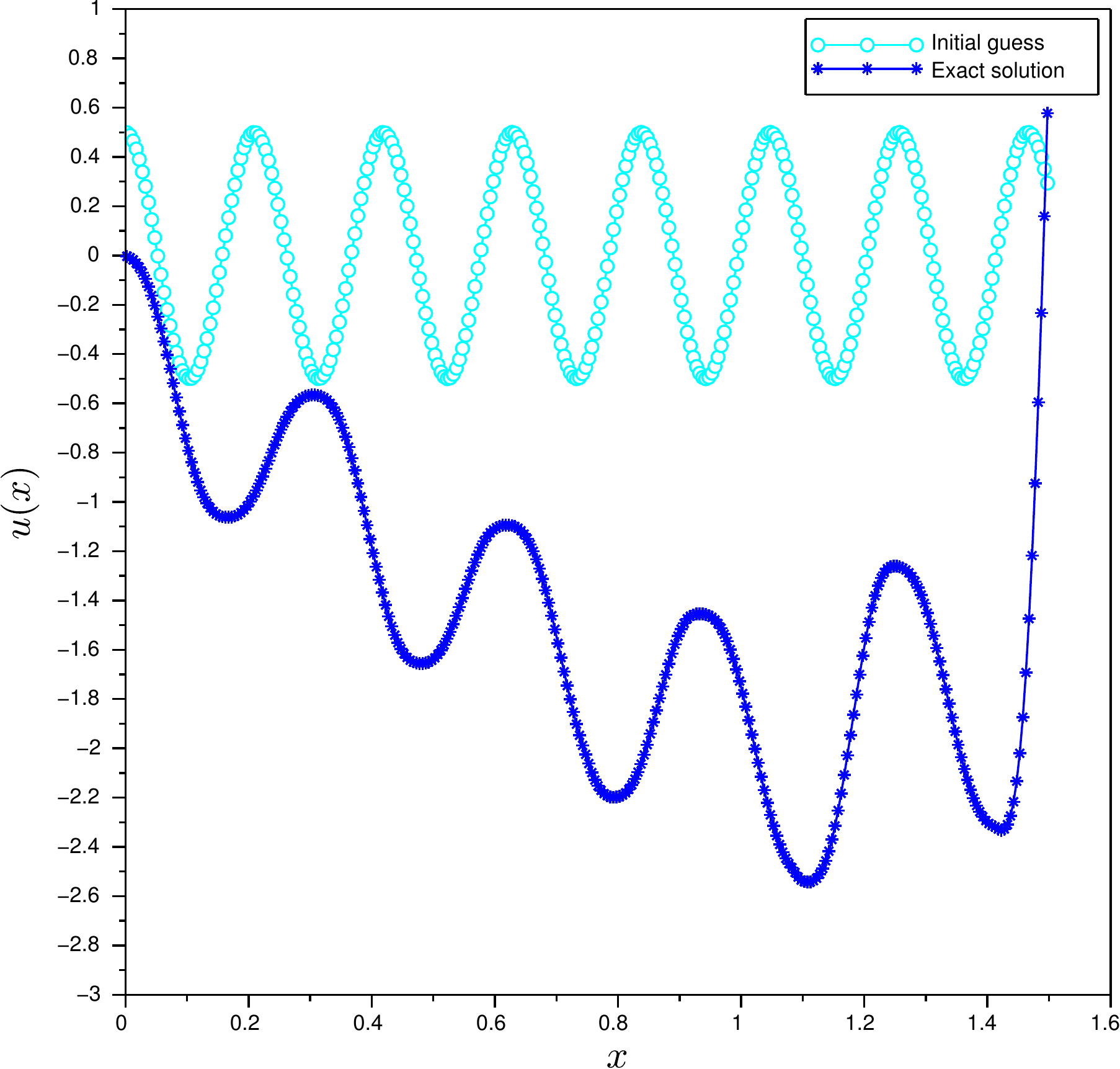}}
  \caption{Permeability field (left), source term (middle), initial
    guess and solution (right).}
  \label{perm}
\end{figure}

We use for the numerical tests the domain $\Omega = (0,3/2)$ with the
boundary conditions $u(0)=0$ and $u(\FK{3/2})=1$, and different values of
$\beta$.
As a first challenging test, we choose the highly variable
permeability field $\lambda$ and the oscillating right hand side
shown in Figure \ref{perm}.
We measure the relative $l^1$ norms of the error obtained at each
Newton iteration as a function of the parallel linear solves $LS_n$
needed in the subdomains per Newton iteration, which is a realistic
measure for the cost of the method. Each Newton iteration requires two
major steps:
\begin{enumerate}
  \item  the evaluation of the fixed point function ${\cal F}$, which
    means solving a non-linear problem in each subdomain. This is done
    using Newton in an inner iteration on each subdomain, and thus
    requires at each inner iteration a linear subdomain solve
    performed in parallel by all subdomains  \FK{(we have used a sparse
    direct solver for the linear subdomain solves in our experiments, but one
    can also use an iterative method if good preconditioners are available)}.
    We denote the maximum
    number of inner iterations needed by the subdomains at the outer
    iteration $j$ by $ls_j^{in}$, and it is the maximum
    which is relevant, because if other subdomains finish earlier,
    they still have to wait for the last one to finish.
  \item the Jacobian matrix needs to be inverted, which we do by
    GMRES, and each GMRES iteration will also need a linear subdomain
    solve per subdomain. We denote the number of linear solves needed
    by GMRES at the outer Newton iteration step $j$ by $ls_j^{G}$.
\end{enumerate}
Hence, the number of linear subdomain solves for the outer Newton
iteration $j$ to complete is $ls_j^{in}+ls_j^G$, and the total number
of linear subdomain solves after $n$ outer Newton iterations is
$LS_n:=\sum_{j=1}^n\left(ls_j^{in}+ls_j^G\right)$.  In all the
numerical tests, we stop the linear GMRES iterations when the relative
residual falls below $10^{-8}$, and the tolerances for the inner and
outer Newton iterations are also set to $10^{-8}$. Adaptive tolerances
could certainly lead to more savings
\cite{el2011guaranteed,ern2013adaptive}, but our purpose here is to
compare the non-linear preconditioners in a fixed setting. The initial
guess we use in all our experiments is shown in Figure \ref{perm} on
the right, together with the solution.

\begin{figure}
  \centering
  \mbox{\includegraphics[width=0.43\textwidth,height=0.43\textwidth]{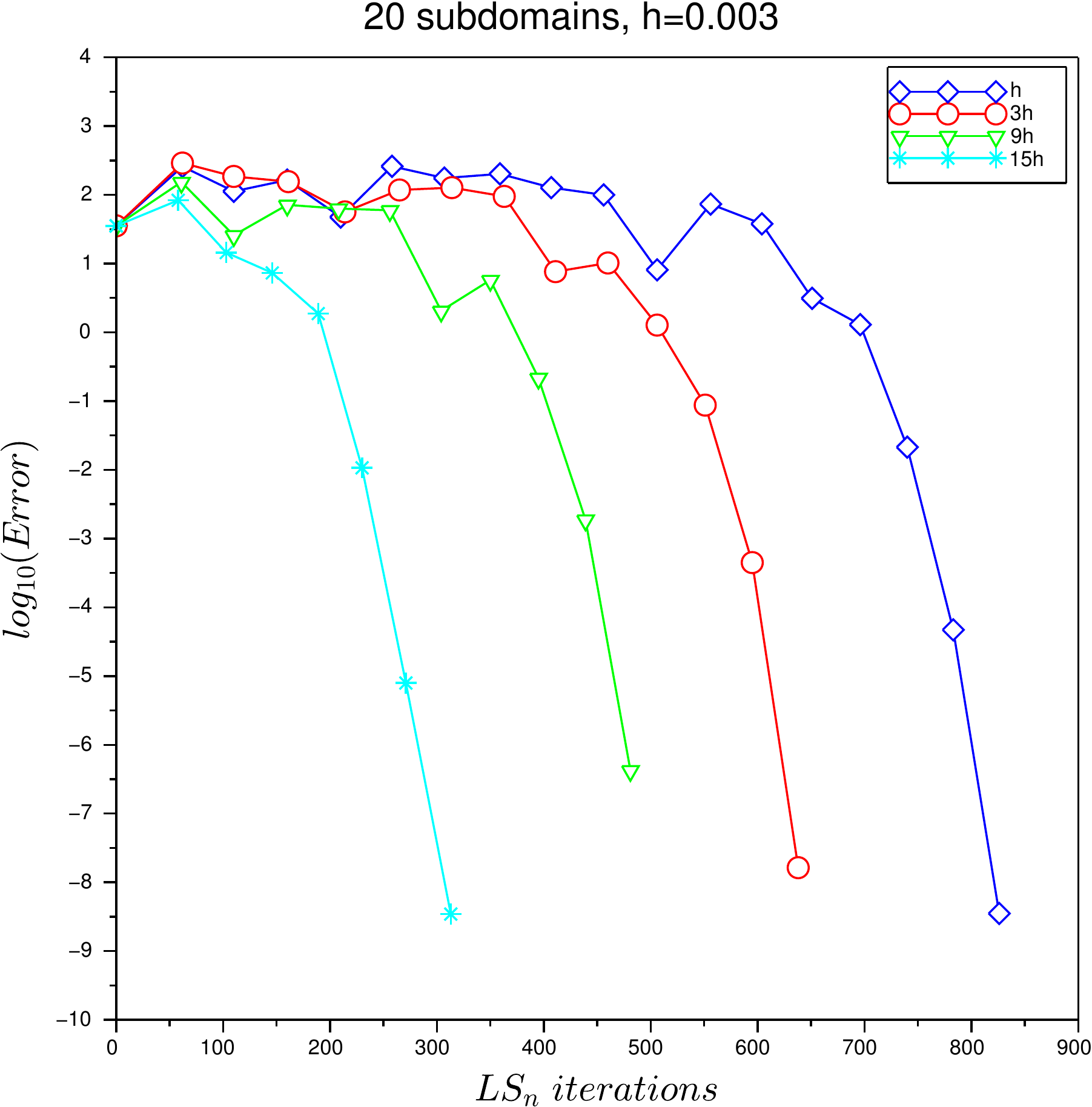}
  \hspace{0.5cm}\includegraphics[width=0.43\textwidth,height=0.43\textwidth]{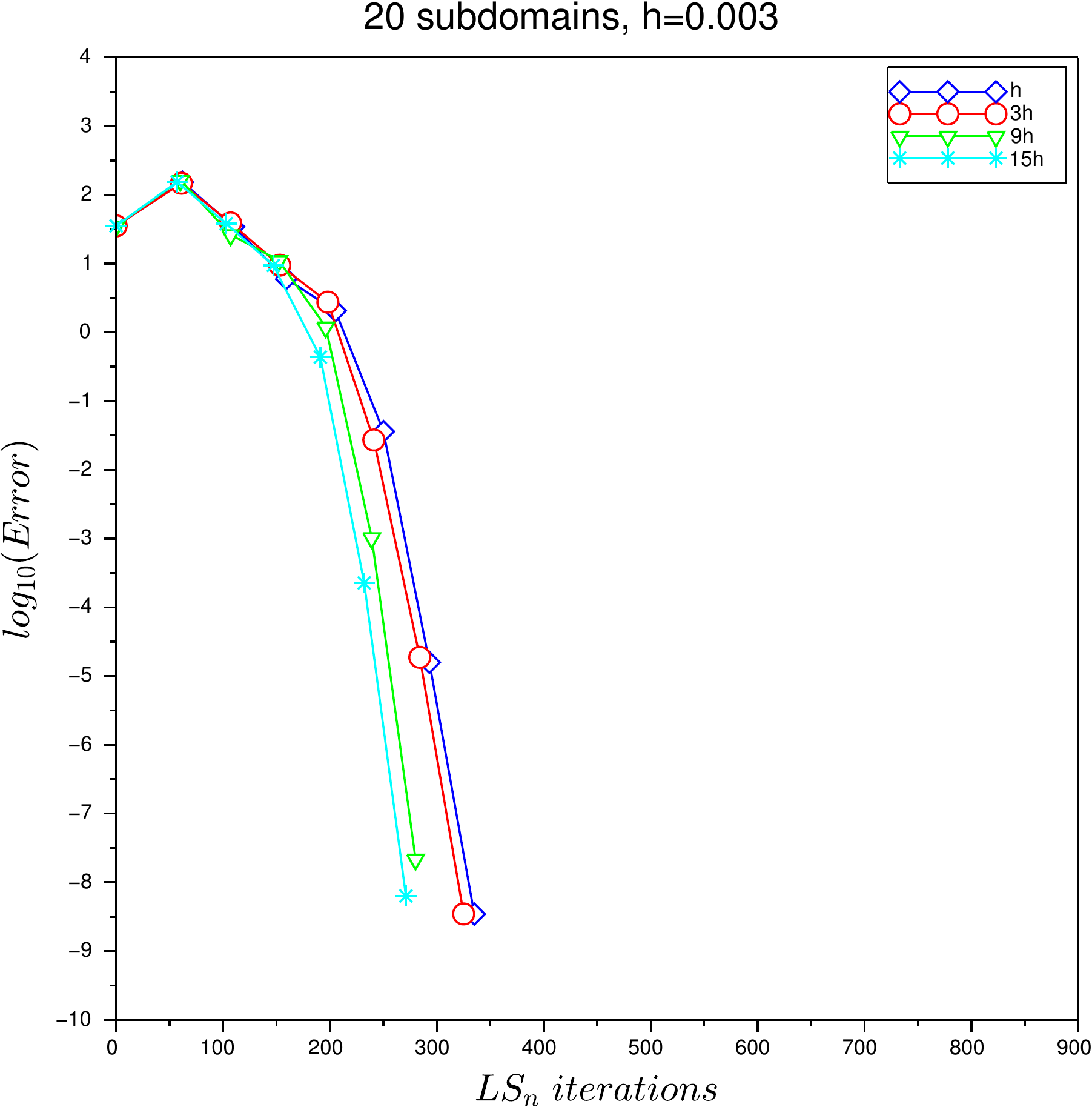}}
  \mbox{
\includegraphics[width=0.43\textwidth,height=0.45\textwidth]{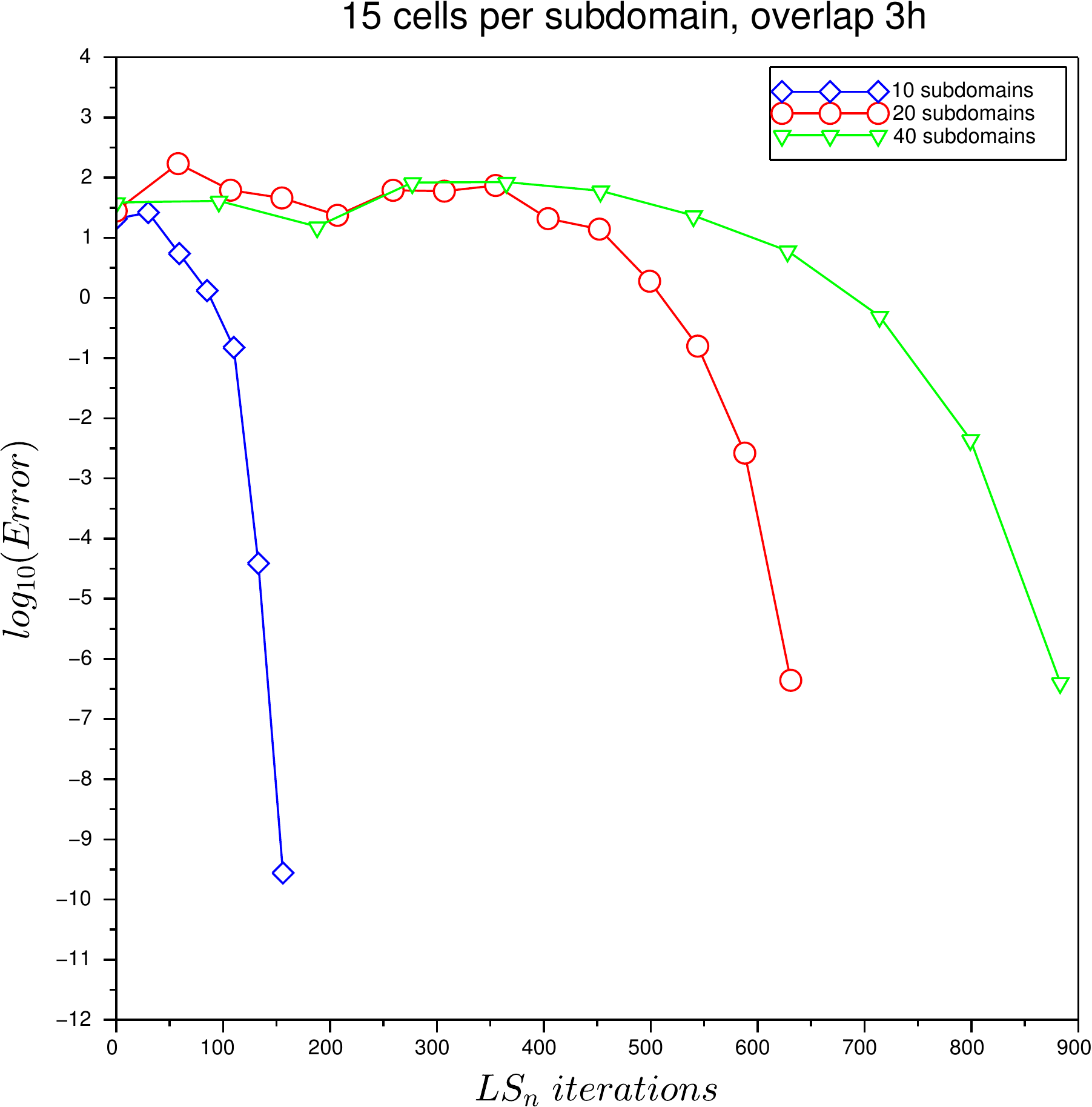}
\hspace{0.5cm}
  \includegraphics[width=0.43\textwidth,height=0.45\textwidth]{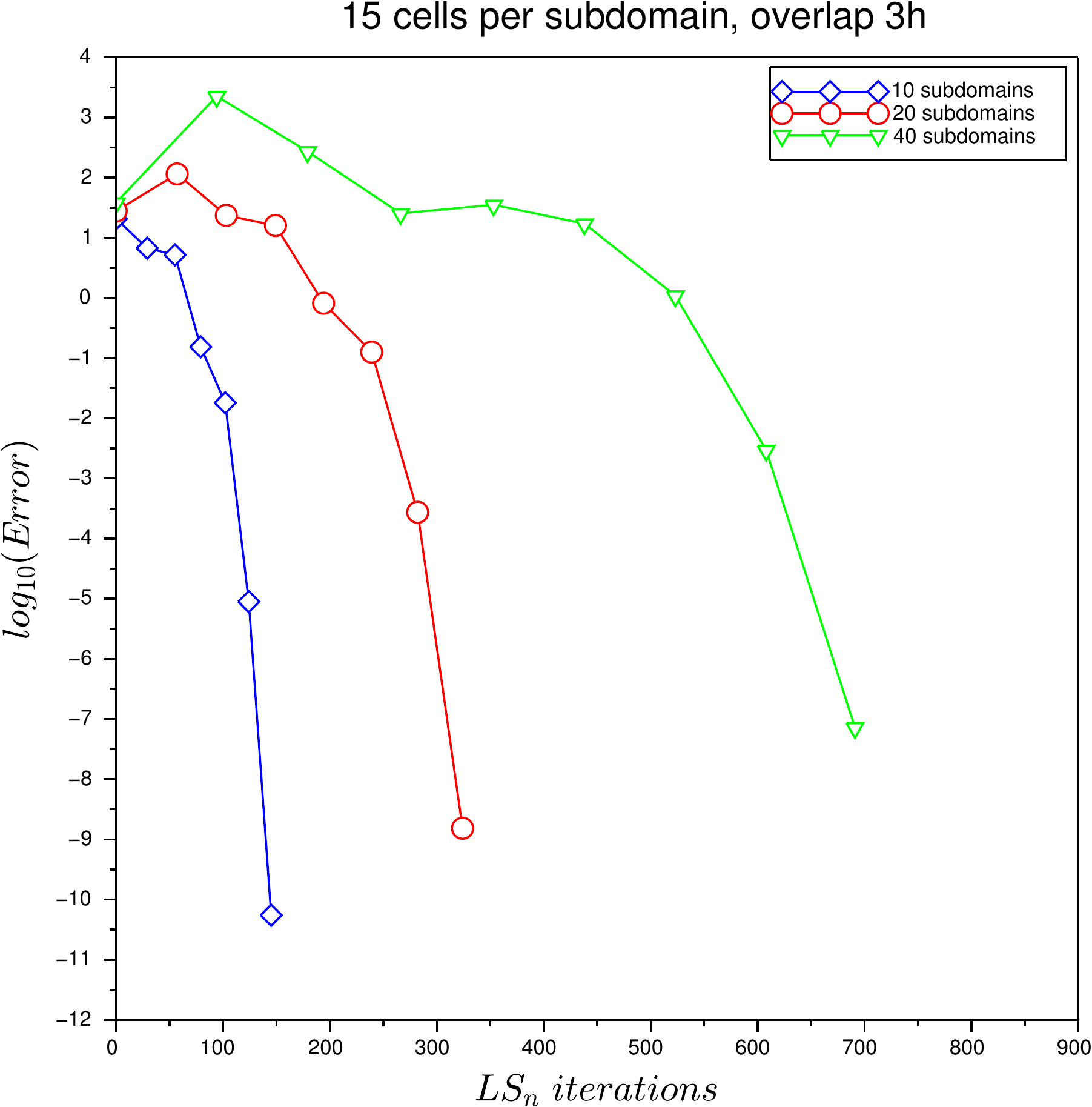}}
  \caption{Error obtained with one-level ASPIN (left) and one-level
    RASPEN (right): in the top row obtained with $20$ subdomains,
    $h=0.003$, and decreasing size of overlap $15h$, $9h$, $3h$, $h$;
    in the bottom row obtained with different number of subdomains
    $10, 20 $ and $40$, overlap $3h$, and a fixed number of cells per
    subdomain.  The Forchheimer problem is defined by the
    permeability, source term, solution and initial guess of
    Figure \ref{perm}.}
  \label{Fig3}
  \centering
  \mbox{\includegraphics[width=0.43\textwidth]{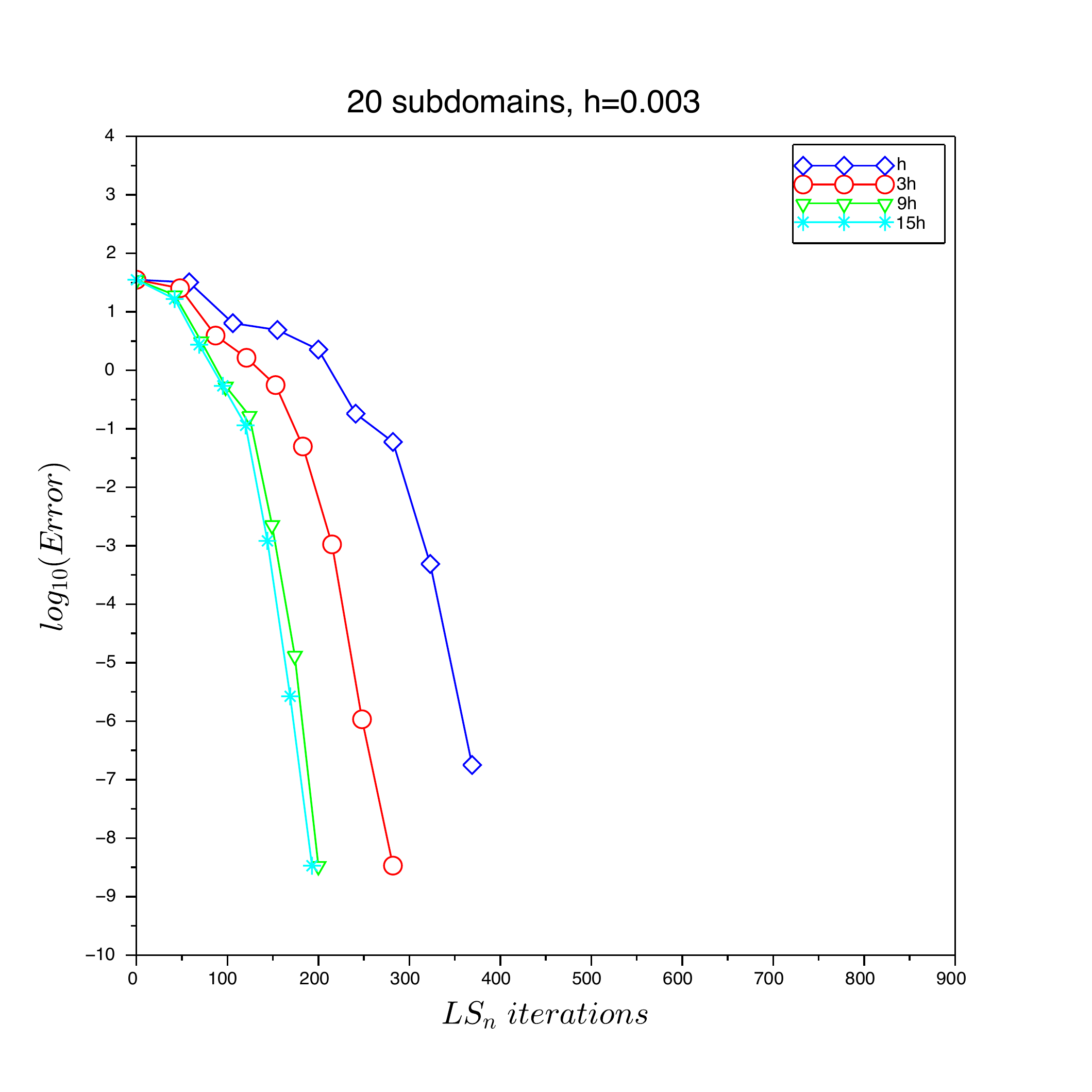}
  \hspace{0.5cm}\includegraphics[width=0.43\textwidth]{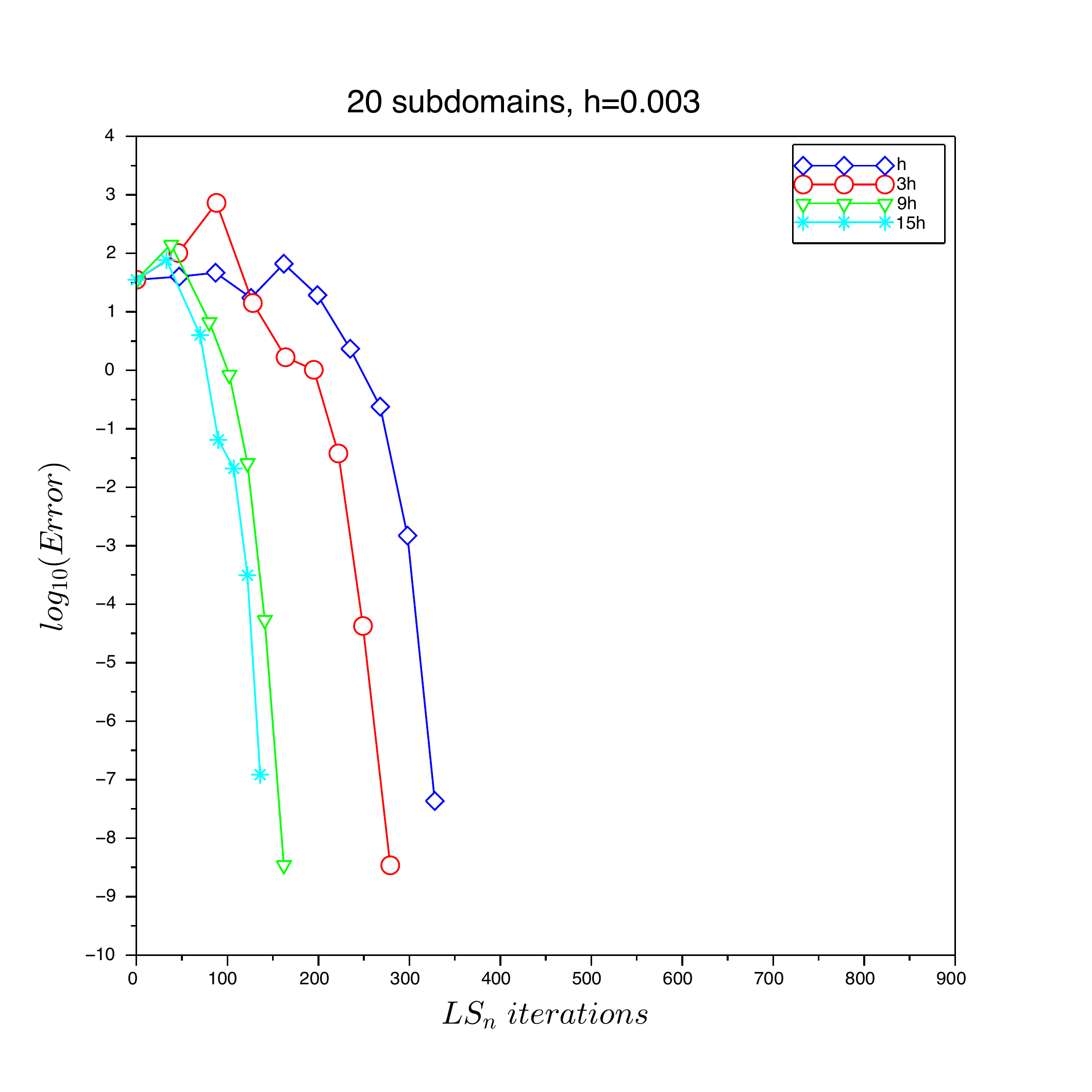}}
  \caption{Error obtained with two-level ASPIN (left) and two-level FAS
    RASPEN (right) obtained with $20$ subdomains, $h=0.003$, and
    decreasing overlap $15h$, $9h$, $3h$, $h$.
The Forchheimer problem is defined by the permeability, source term, solution and initial guess
of Figure \ref{perm}.}
 \label{Fig6}
\end{figure}

We show in Figure \ref{Fig3}
how the convergence depends on the overlap and the number of
subdomains for one level ASPIN and RASPEN with Forchheimer model
parameter $\beta=1$.
In the top row on the left of Figure \ref{Fig3},
we see that for ASPIN the number of linear iterations increases much
more rapidly when decreasing the overlap than for RASPEN on the right
for a fixed mesh size $h=0.003$ and number of subdmains equal $20$. In
the bottom row of Figure \ref{Fig3}, we see that the convergence of
both one level ASPIN and RASPEN depends on the number of subdomains,
but RASPEN seems to be less sensitive than ASPIN.


\subsubsection{Two level variants}

In Figure \ref{Fig6}, we show the dependence of two-level ASPIN and
two-level FAS-RASPEN on a decreasing size of the overlap, as we did
for the one-level variants in the top row of Figure \ref{Fig3}.
We see that the addition of the coarse level improves the performance,
for RASPEN when the overlap is large, and in all cases for ASPIN.

\begin{figure}
  \centering
  \mbox{\includegraphics[width=0.45\textwidth]{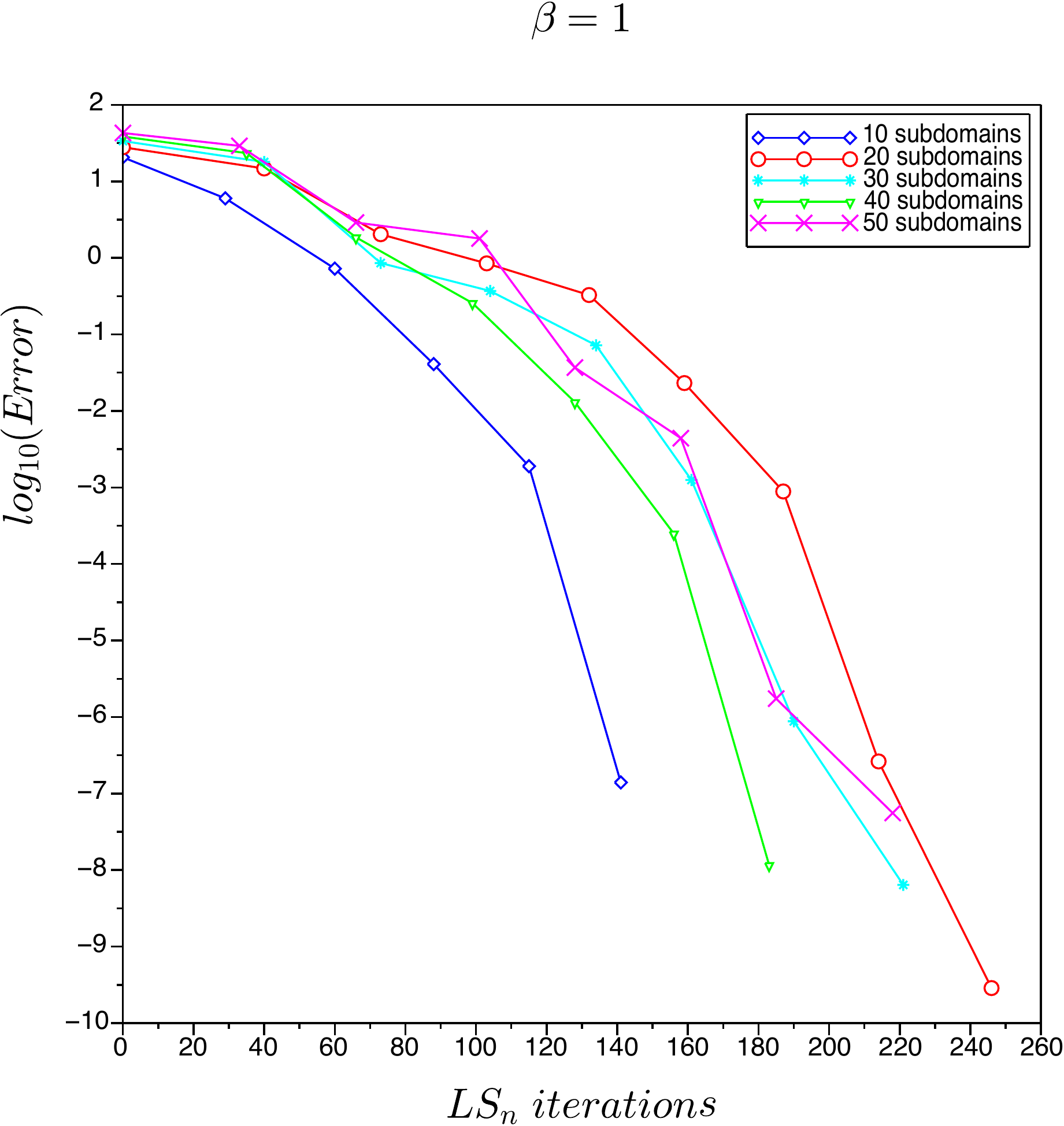}
\includegraphics[width=0.45\textwidth]{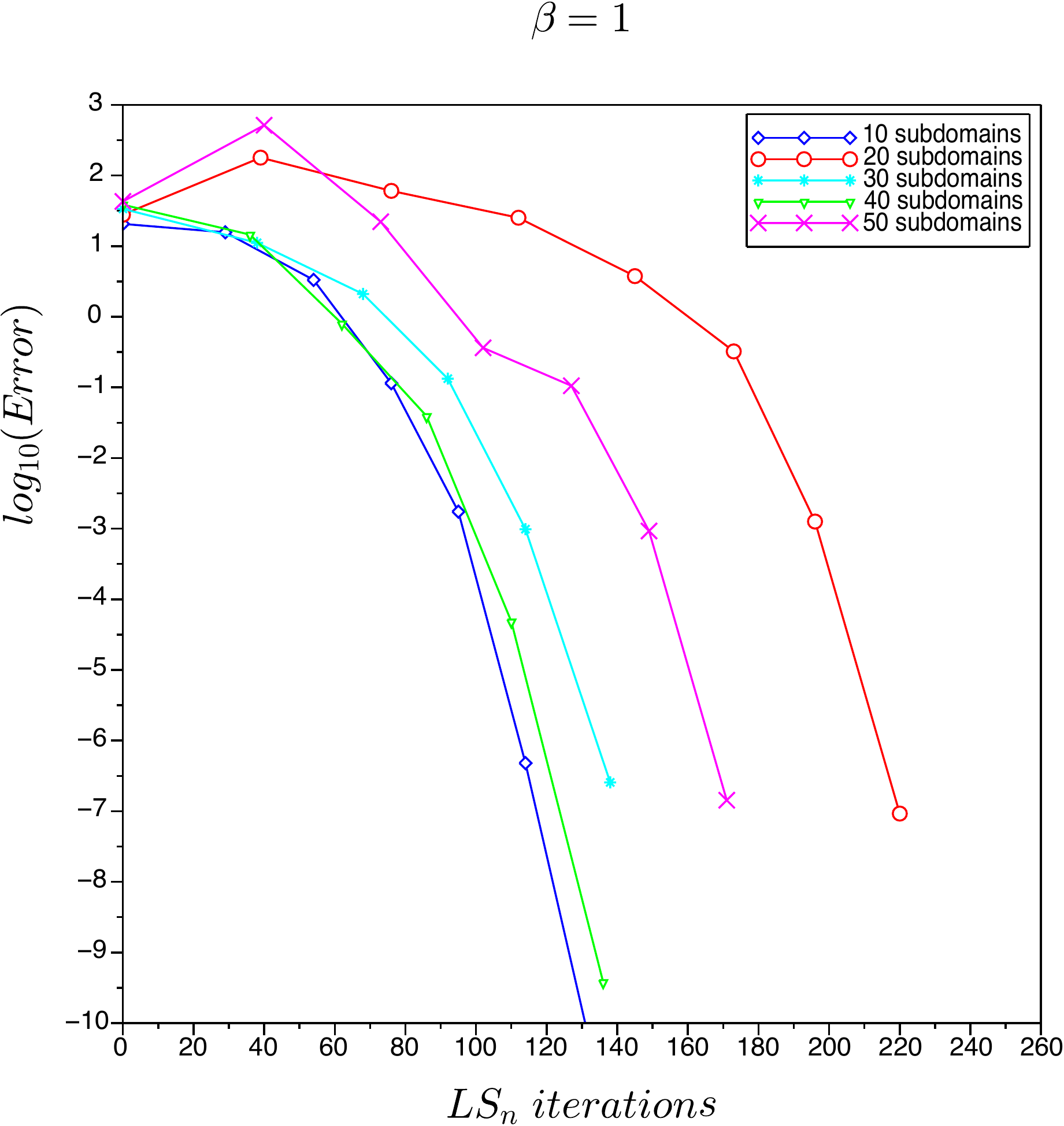}}
 \mbox{ \includegraphics[width=0.45\textwidth]{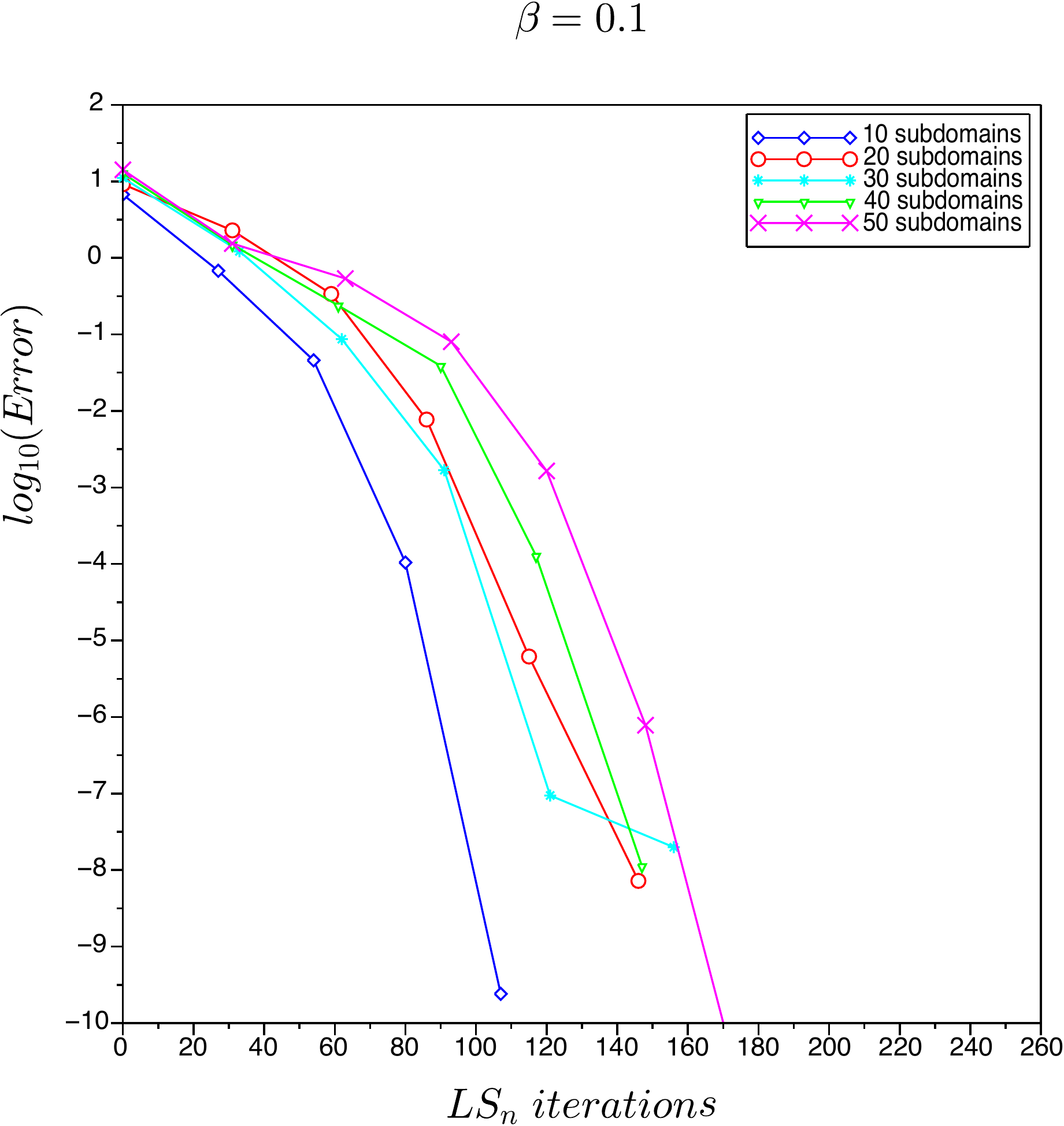}
 \includegraphics[width=0.45\textwidth]{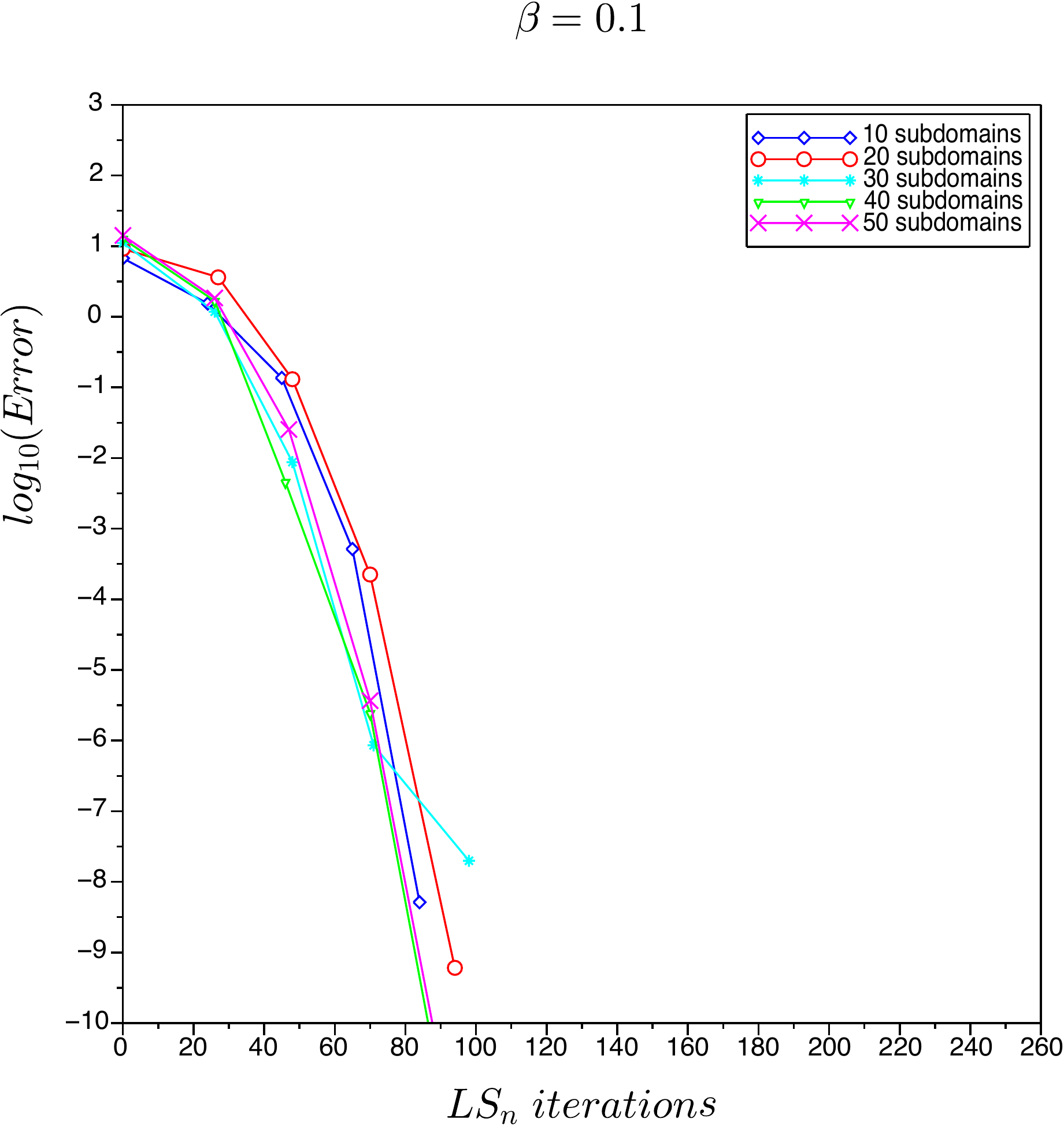}}
  \mbox{\includegraphics[width=0.45\textwidth]{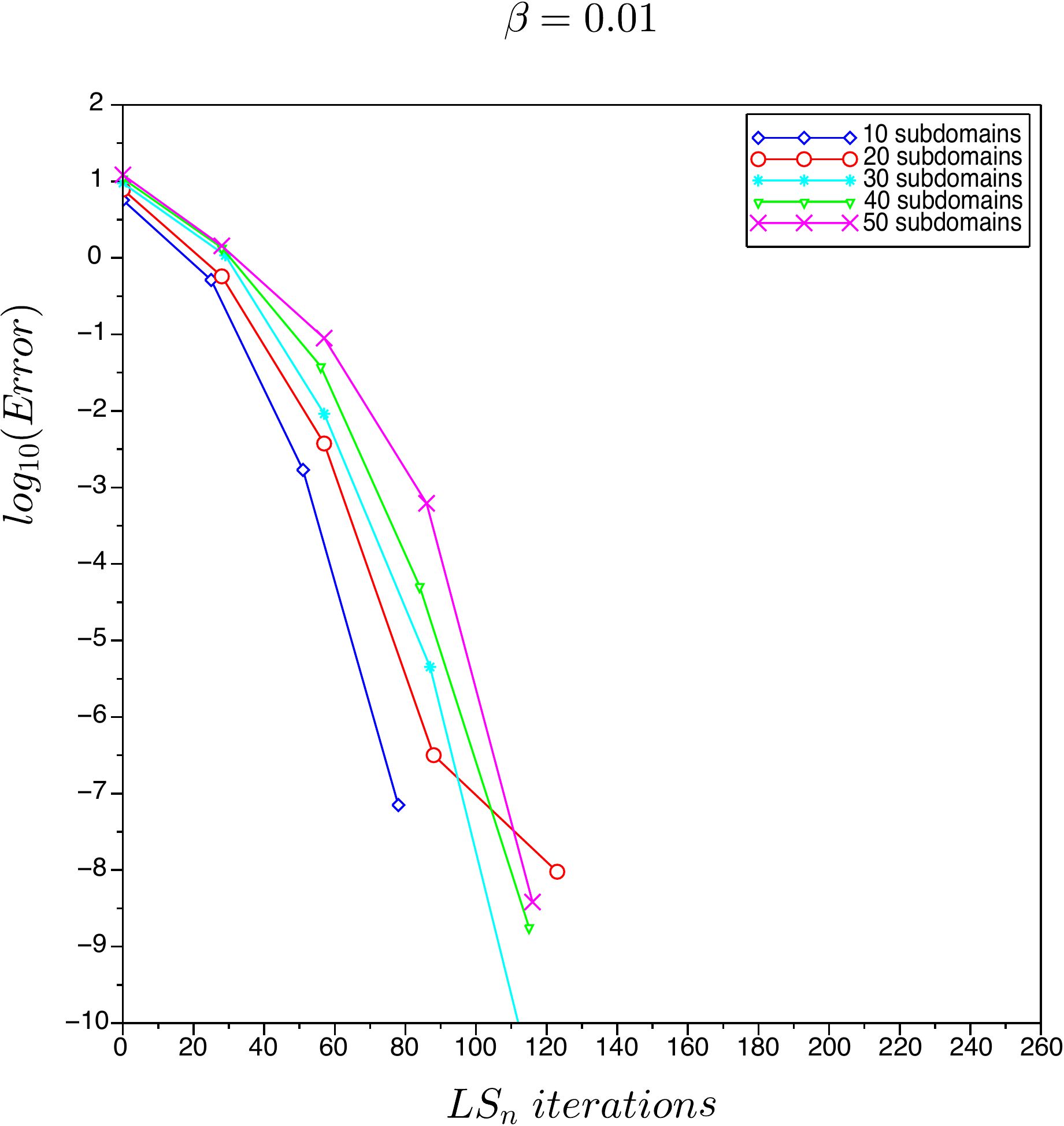}
\includegraphics[width=0.45\textwidth]{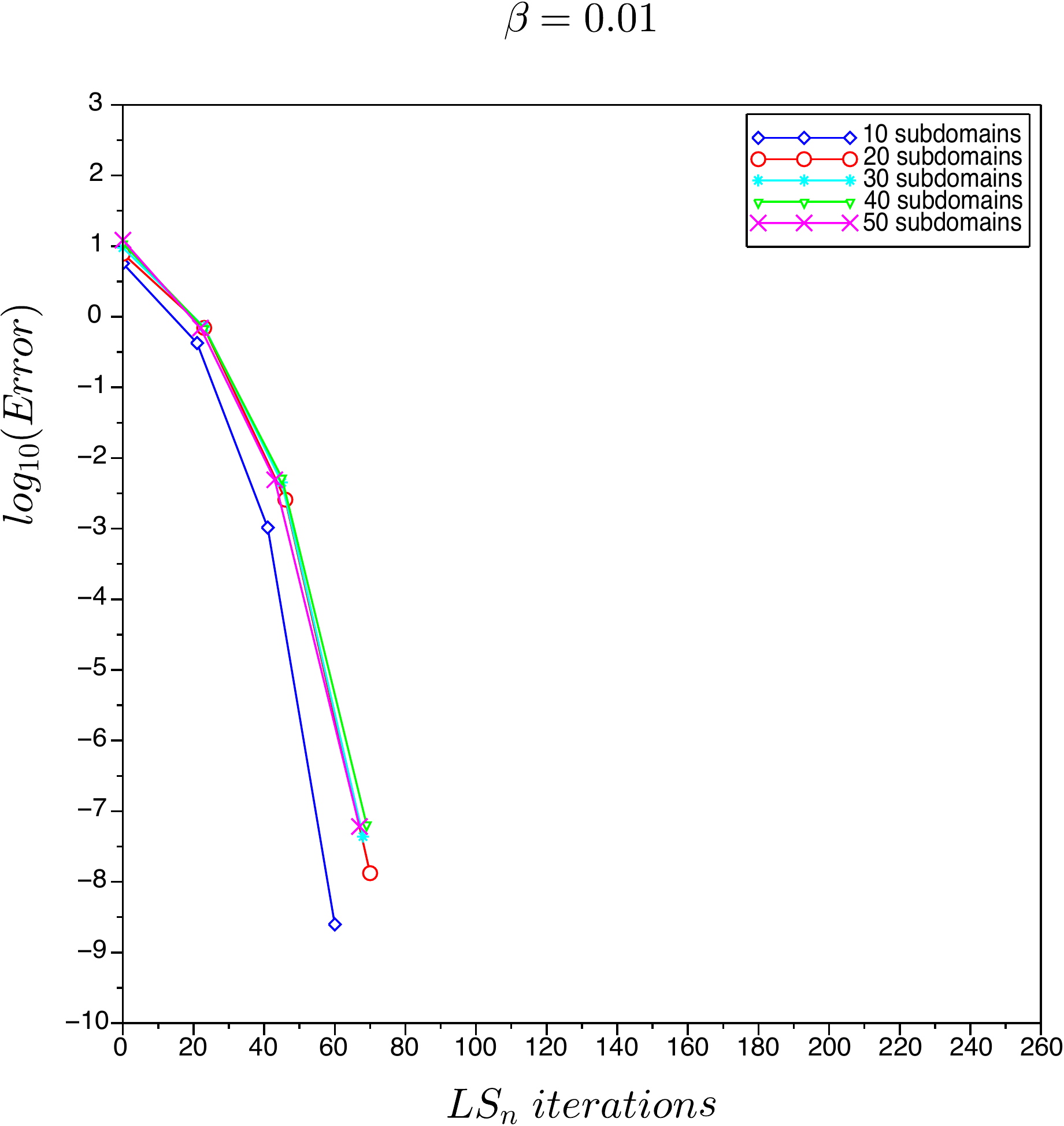}}
  \caption{Error obtained with two-level ASPIN (left) and two-level
    FAS RASPEN (right) and different number of subdomains $10, 20, 30,
    40, 50$.  From top to bottom with decreasing Forchheimer parameter
    $\beta=1,0.1,0.01$. The Forchheimer problem is defined by the
    permeability, source term, solution and initial guess of Figure
    \ref{perm}. }
  \label{Fig6-2}
\end{figure}

\begin{table}\tabcolsep0.5em
\begin{tabular}{|l|l|l|l|l|l|l|l|l|}
\hline
\multicolumn{7}{|l|}{ASPIN}                                                                           \\ \hline
 \footnotesize Number of subdomains & \multicolumn{2}{l|}{10} & \multicolumn{2}{l|}{20} & \multicolumn{2}{l|}{40} \\ \hline
   \theadfont\diagbox[width=11em]{ \footnotesize Overlap size}{ \footnotesize type of iteration}  &    \footnotesize PIN iter.  &    \footnotesize $LS_n$ iter. &      \footnotesize PIN iter.  &        \footnotesize  $LS_n$ iter.    &   \footnotesize PIN iter.  &        \footnotesize  $LS_n$ iter.   \\ \hline
  h &      8     &      184     &      15     &     663      &      -     &           -       \\ \hline
 3h&      7     &       156    &     14      &     631      &        11   &    883            \\ \hline
 5h&      6     &     130      &    11       &     479      &        10   &      744             \\ \hline
\multicolumn{7}{|l|}{RASPEN}                                                                           \\ \hline
 \footnotesize Number of subdomains & \multicolumn{2}{l|}{10} & \multicolumn{2}{l|}{20} & \multicolumn{2}{l|}{40} \\ \hline
   \theadfont\diagbox[width=11em]{ \footnotesize Overlap size}{ \footnotesize type of iteration}  &  \footnotesize PEN iter.  &        \footnotesize $LS_n$ iter. &     \footnotesize PEN iter.  &      \footnotesize $LS_n$ iter.    &   \footnotesize PEN iter.  &      \footnotesize $LS_n$ iter.   \\ \hline
  h &    7       &     150      &     9      &    369       &     9      &       701           \\ \hline
 3h&     7      &     145      &     8      &    324       &    9       &     691           \\ \hline
 5h&      6     &     126      &     7      &     274     &     9       &       659            \\ \hline
\multicolumn{7}{|l|}{ \footnotesize Two-level  ASPIN}                                                                           \\ \hline
 \footnotesize Number of subdomains & \multicolumn{2}{l|}{10} & \multicolumn{2}{l|}{20} & \multicolumn{2}{l|}{40} \\ \hline
\theadfont\diagbox[width=11em]{ \footnotesize Overlap size}{ \footnotesize type of iteration}  &   \footnotesize PIN iter.  &      \footnotesize $LS_n$ iter. &     \footnotesize PIN iter.  &      \footnotesize $LS_n$ iter.    &   \footnotesize PIN iter.  &      \footnotesize $LS_n$ iter.   \\ \hline  h &    7       &     184      &    9       &     316      &    8       &      285            \\ \hline
3h &     6      &      141     &      9     &    246       &    7       &        183           \\ \hline
5h &      6     &   135        &      8     &       199    &     7      &        164          \\ \hline
 \multicolumn{7}{|l|}{ \footnotesize Two-level FAS-RASPEN}                                                                           \\ \hline
 \footnotesize Number of subdomains & \multicolumn{2}{l|}{10} & \multicolumn{2}{l|}{20} & \multicolumn{2}{l|}{40} \\ \hline
   \theadfont\diagbox[width=11em]{ \footnotesize Overlap size}{ \footnotesize type of iteration}  &   \footnotesize PEN iter.  &      \footnotesize $LS_n$ iter. &     \footnotesize PEN iter.  &      \footnotesize $LS_n$ iter.    &   \footnotesize PEN iter.  &      \footnotesize $LS_n$ iter.   \\ \hline  h &      7     &      134     &     9      &    272       &     8      &        258          \\ \hline
 3h&     7      &     133      &      8     &    220       &      6     &     136           \\ \hline
 5h&     6      &     112      &      8     &     211    &       6    &         116          \\ \hline
\end{tabular}
\caption{Comparison in terms of non-linear and linear iterations of the
different algorithms for the Forchheimer problem defined by
the permeability, source term, solution and initial guess
of Figure \ref{perm}. \label{table:questions:1}}
\end{table}

\begin{figure}[!ht]
  \centering
  \mbox{\includegraphics[width=0.48\textwidth]{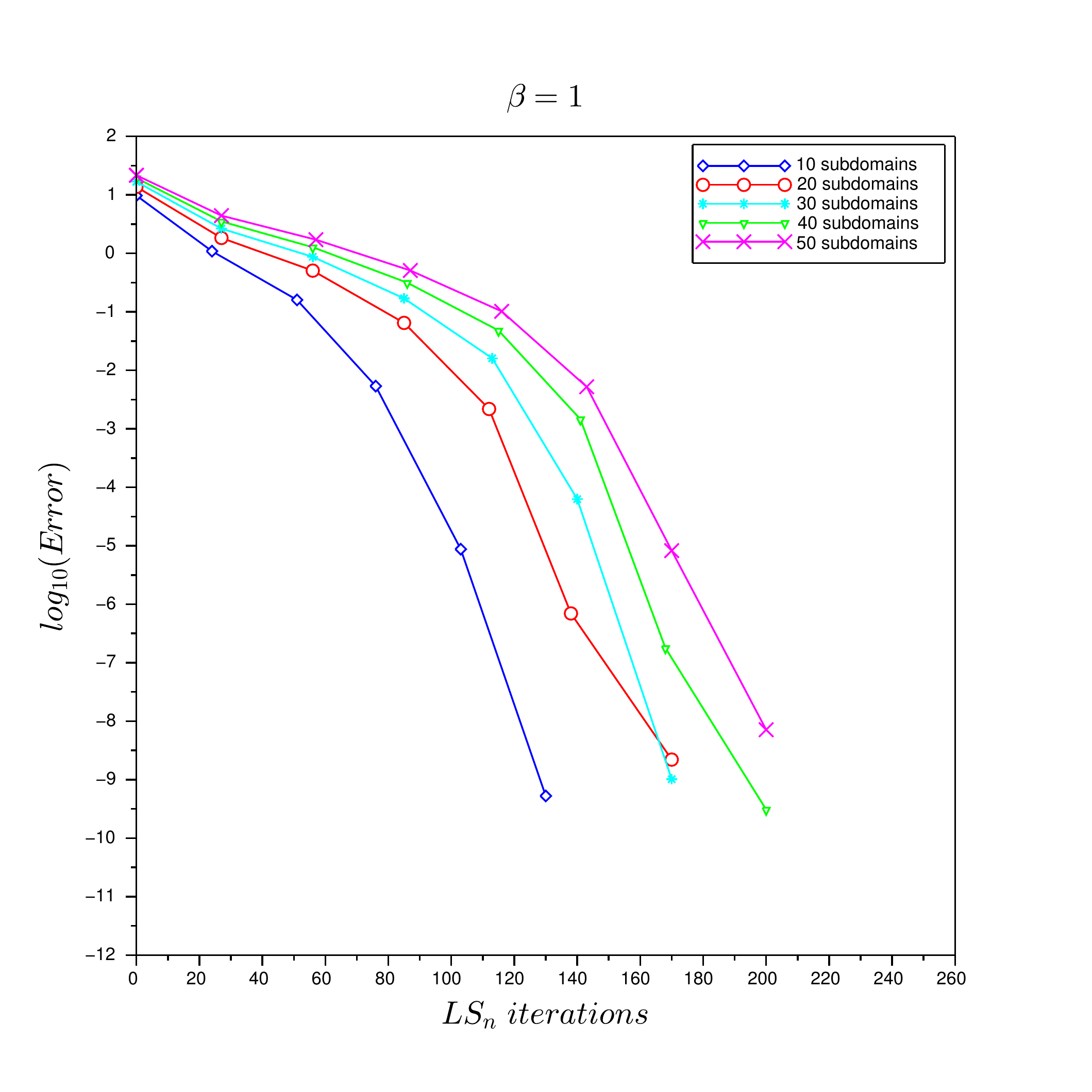}
\includegraphics[width=0.48\textwidth]{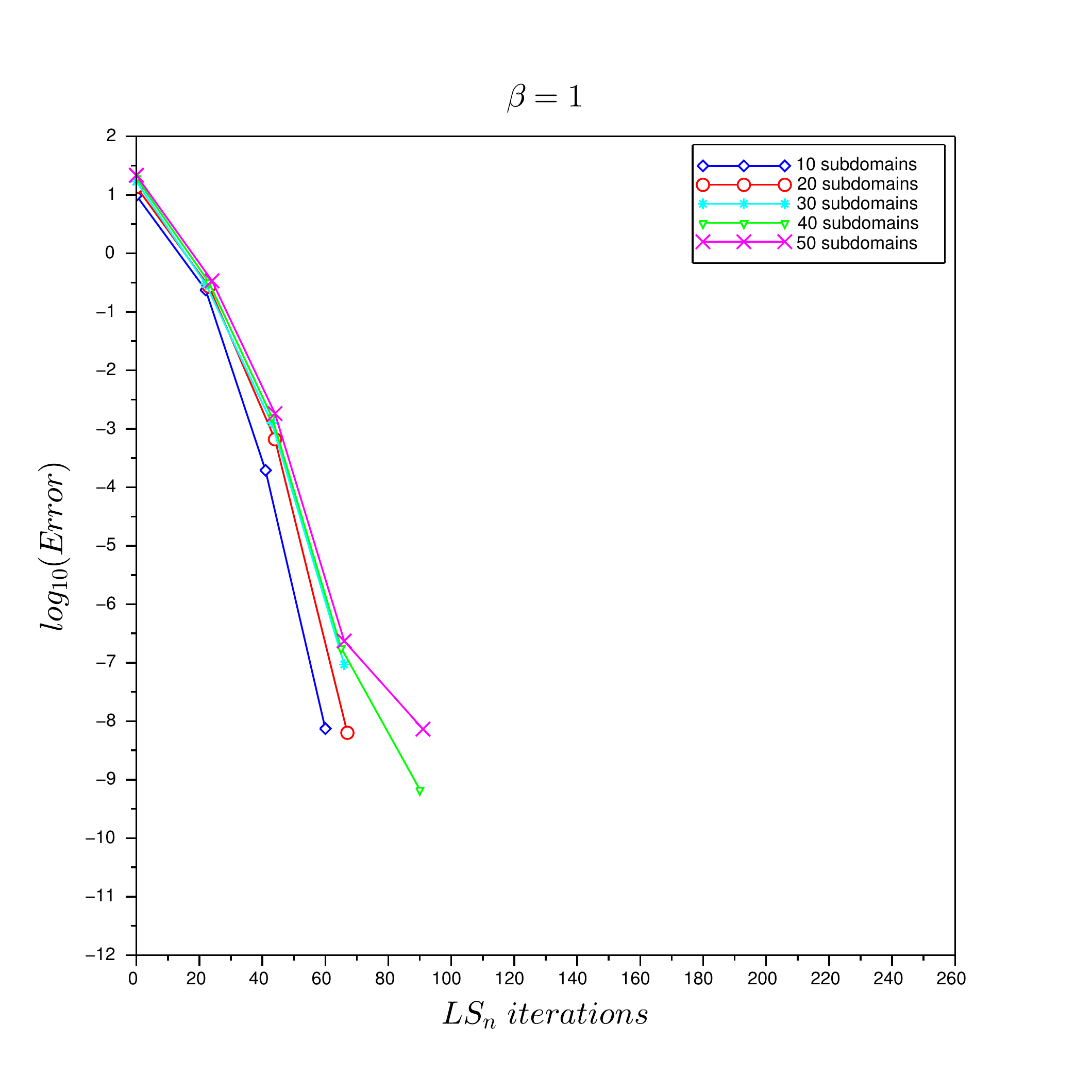}}\\*[-15pt]
  \caption{Error obtained with two-level ASPIN (left) and two-level
    FAS RASPEN (right) with overlap $3h$, and different number of subdomains $10, 20, 30,
    40, 50$ for the smooth Forchheimer example. \\*[-30pt]}
  \label{Figb=1}
\end{figure}
\begin{table}[!hb]
{\small\tabcolsep0.5em
\begin{tabular}{|l|l|l|l|l|l|l|l|l|}
\hline
\multicolumn{7}{|l|}{ASPIN}                                                                           \\ \hline
 \footnotesize Number of subdomains & \multicolumn{2}{l|}{10} & \multicolumn{2}{l|}{20} & \multicolumn{2}{l|}{40} \\ \hline
   \theadfont\diagbox[width=11em]{ \footnotesize Overlap size}{ \footnotesize type of iteration}  &    \footnotesize PIN iter.  &    \footnotesize $LS_n$ iter. &      \footnotesize PIN iter.  &        \footnotesize  $LS_n$ iter.    &   \footnotesize PIN iter.  &        \footnotesize  $LS_n$ iter.   \\ \hline
  h &      5     &      118     &      5     &     228      &      6     &           520       \\ \hline
 3h&      5     &      118    &     5      &     227      &        6   &       516            \\ \hline
 5h&      5     &      117      &    5       &     222      &        6   &      480             \\ \hline
\multicolumn{7}{|l|}{RASPEN}                                                                           \\ \hline
 \footnotesize Number of subdomains & \multicolumn{2}{l|}{10} & \multicolumn{2}{l|}{20} & \multicolumn{2}{l|}{40} \\ \hline
   \theadfont\diagbox[width=11em]{ \footnotesize Overlap size}{ \footnotesize type of iteration}  &  \footnotesize PEN iter.  &        \footnotesize $LS_n$ iter. &     \footnotesize PEN iter.  &      \footnotesize $LS_n$ iter.    &   \footnotesize PEN iter.  &      \footnotesize $LS_n$ iter.   \\ \hline
  h &    4      &     92      &     4      &    172       &     4      &      340           \\ \hline
 3h&     4      &    87      &     4      &    172        &    4       &    331           \\ \hline
 5h&      4     &     88      &     4      &     168     &     4       &       313           \\ \hline
\multicolumn{7}{|l|}{ \footnotesize Two level  ASPIN}                                                                           \\ \hline
 \footnotesize Number of subdomains & \multicolumn{2}{l|}{10} & \multicolumn{2}{l|}{20} & \multicolumn{2}{l|}{40} \\ \hline
\theadfont\diagbox[width=11em]{ \footnotesize Overlap size}{ \footnotesize type of iteration}  &   \footnotesize PIN iter.  &      \footnotesize $LS_n$ iter. &     \footnotesize PIN iter.  &      \footnotesize $LS_n$ iter.    &   \footnotesize PIN iter.  &      \footnotesize $LS_n$ iter.   \\ \hline
  h &    5       &     140      &    5       &     240      &    5       &      280           \\ \hline
3h &     5      &      130     &      6     &    170       &    6       &        200          \\ \hline
5h &      5     &   115        &      7     &       149    &     6      &        147         \\ \hline
 \multicolumn{7}{|l|}{ \footnotesize Two level FAS RASPEN}                                                                           \\ \hline
 \footnotesize Number of subdomains & \multicolumn{2}{l|}{10} & \multicolumn{2}{l|}{20} & \multicolumn{2}{l|}{40} \\ \hline
   \theadfont\diagbox[width=11em]{ \footnotesize Overlap size}{ \footnotesize type of iteration}  &   \footnotesize PEN iter.  &      \footnotesize $LS_n$ iter. &     \footnotesize PEN iter.  &      \footnotesize $LS_n$ iter.    &   \footnotesize PEN iter.  &      \footnotesize $LS_n$ iter.   \\ \hline
  h &      4     &      77     &     3      &    87       &     4      &        131          \\ \hline
 3h&     3      &     60      &      3     &    67      &      4    &     90           \\ \hline
 5h&     3      &     55      &      3     &     57    &       3    &         57          \\ \hline
\end{tabular}
}
\caption{Comparison in terms of non-linear and linear iterations of the
different algorithms for the smooth Forchheimer example.}
\label{table:questions:2}
\end{table}




\begin{table}
\centering
{\small\tabcolsep0.5em
\begin{tabular}{|c|c|c|c|c|c|c|c|c|c|}
  \hline
  Number of &  & \multicolumn{4}{c|}{1-Level} & \multicolumn{4}{c|}{2-Level}\\
  \cline{3-10}
subdomains & $n$  &  $ls_n^G$  & $ls_n^{in}$ & $ls_n^{\min}$ &
$LS_n$ & $ls_n^G$  & $ls_n^{in}$ & $ls_n^{\min}$ &
$LS_n$
\\\hline
$10$ & 1 & 19 (20) & 4 (4) & 3 (3) &         & 15 (20) & 7 (4) & 3 (3) &          \\
     & 2 & 19 (20) & 3 (6) & 3( 3) & 87 (118)& 16 (21) & 3 (6) & 2 (3) & 60 (130) \\
     & 3 & 19 (20) & 2 (4) & 2 (2) &         & 17 (22) & 2 (3) & 1 (2) &          \\
     & 4 & 19 (20) & 2 (2) & 1 (2) &         & - (24) & - (3) & - (1) &           \\
     & 5 &  -  (21) & - (1) & - (1) &        & - (25) & - (2) & - (1) &           \\
\hline
$20$ & 1 & 40 (41) & 5 (5) & 3 (3) &           & 15 (22) & 8 (5) & 3 (3) &         \\
     & 2 & 40 (41) & 3 (7) & 2 (2) & 172 (227) & 18 (23) & 3 (6) & 2 (3) & 67 (170)\\
     & 3 & 40 (41) & 2 (5) & 1 (2) &           & 21 (24) & 2 (5) & 1 (2) &         \\
     & 4 & 40 (41) & 2 (3) & 1 (1) &           & - (24) & - (3) & - (1) &          \\
     & 5 & -  (41) & - (2) & - (1) &           & - (24) & - (2) & - (1) &          \\
     & 6 & -  (-)  & - (-) & - (-) &           & - (31) & - (1) & - (1) &          \\

\hline
$40$ & 1 & 78 (80) & 5 (5) & 3 (3) &          & 14 (22) & 9 (5) & 3 (3) &          \\
     & 2 & 81 (81) & 3 (6) & 2 (2) & 331 (516)& 17 (22) & 3 (7) & 1 (2) & 90 (200) \\
     & 3 & 79 (82) & 2 (6) & 1 (2) &          & 20 (24) & 2 (6) & 1 (2) &          \\
     & 4 & 81 (82) & 2 (5) & 1 (1) &          & 24 (24) & 1(5) & 0 (1) &           \\
     & 5 & -  (82) & - (3) & - (1) &          & - (23) & - (3) & - (1) &           \\
     & 6 & -  (82) & - (2) & - (1) &          & - (25) & - (2) & - (1) &           \\
     & 7 & -  (-)  & - (-) & - (-) &          & - (31) & - (1) & - (0) &           \\

\hline
\end{tabular}
}
\caption{Numerical results with one- and two-level RASPEN and ASPIN for the
  1D non-linear smooth Forchheimer problem. `-' indicates that the method has converged.}
\label{table:results1levFroDisc}
\end{table}

%
%
%
%
%
%

In Figure \ref{Fig6-2},
we present a study of the influence of the number of subdomains on the
convergence for two-level ASPIN and two-level FAS-RASPEN with
different values of the Forchheimer parameter $\beta=1, 0.1, 0.01$
which governs the non-linearity of the model (the model becomes linear
for $\beta=0$). An interesting observation is that for $\beta=1$, the convergence
of both two-level ASPIN and two-level
FAS-RASPEN depends on the number of subdomains in an irregular fashion:
increasing the number of subdomains sometimes increases iteration
counts, and then decreases them again. We will study this effect further below, but note
already from Figure \ref{Fig6-2} that this dependence disappears for two-level FAS-RAPSEN as the
the nonlinearity diminishes (i.e., as $\beta$ decreases), and is weakened for two-level ASPIN.

We finally show in Table \ref{table:questions:1}
the number of outer Newton iterations (PIN iter for ASPIN and PEN iter
for RASPEN) and the total number of linear iterations ($LS_n$ iter) for
various numbers of subdomains and various overlap sizes obtained with ASPIN,
RASPEN, two-level ASPIN and two-level FAS-RASPEN. We see that the
coarse grid considerably improves the convergence of both RASPEN and
ASPIN. Also, in all cases, RASPEN needs substantially
fewer linear iterations than ASPIN.

We now return to the irregular number of iterations observed in
  Figure \ref{Fig6-2} for the Forchheimer parameter $\beta=1$, i.e when
  the non-linearity is strong. We claim that this irregular dependence is
  due to strong variations in the initial guesses used by RASPEN and
  ASPIN at subdomain interfaces, which is in turn caused by the highly
  variable contrast and oscillating source term we used, leading to an oscillatory
  solution, see Figure \ref{perm}. In other words, we expect the irregularity
  to disappear when the solution is non-oscillatory. To test this, we now present
  numerical results with the less variable permeability function
  $\lambda(x)=\cos(x)$ and source term $f(x)=\cos(x)$ as well, which leads
  to a smooth solution. Starting with a zero initial guess, we show in
  Figure \ref{Figb=1} the results obtained for Forchheimer parameter
  $\beta=1$, corresponding to the first row of Figure \ref{Fig6-2}.

We clearly see that the irregular behavior has now disappeared for
  both two-level ASPIN and RASPEN, but two-level ASPIN still shows
  some dependence of the iteration numbers as the number of subdomains
  increases. We show in Table \ref{table:questions:2}
the complete results for this smoother example, and we see that
  the irregular convergence behavior of the two-level methods is no longer
   present. We finally give in Table \ref{table:results1levFroDisc}
a detailed account of the linear subdomain solves needed for each outer
  Newton iteration $n$ for the case of an overlap of $3h$.
  There, we use the format
  $it_{\mathrm{RASPEN}}(it_{\mathrm{ASPIN}})$, where
  $it_{\mathrm{RASPEN}}$ is the iteration count for RASPEN and
  $it_{\mathrm{ASPIN}}$ is the iteration count for ASPIN. We show in
  the first column the linear subdomain solves $ls_n^G$ required for
  the inversion of the Jacobian matrix using GMRES, see item 2 in
  Subsection \ref{subseconelevelnum}, and in the next column the
  maximum number of iterations $ls_n^{in}$ needed to evaluate the
  nonlinear fixed point function ${\cal F}$, see item 1 in Subsection
  \ref{subseconelevelnum}. In the next column, we show for
  completeness also the smallest number of inner iterations
  $ls_n^{\min}$ any of the subdomains needed, to illustrate how
  balanced the work is in this example. The last column then contains the
  total number of linear iterations $LS_n$, see Subsection
  \ref{subseconelevelnum}. These results show that the main gain of
  RASPEN is a reduced number of Newton iterations, i.e. it is a better
  non-linear preconditioner than ASPIN, and also a reduced number of
  inner iterations for the non-linear subdomain solves, i.e. the
  preconditioner is less expensive. This leads to the substantial
  savings observed in the last columns, and in Table \ref{table:questions:2}.

\subsection{A non-linear Poisson problem}

\begin{figure}
  \centering
  \includegraphics[width=0.6\textwidth]{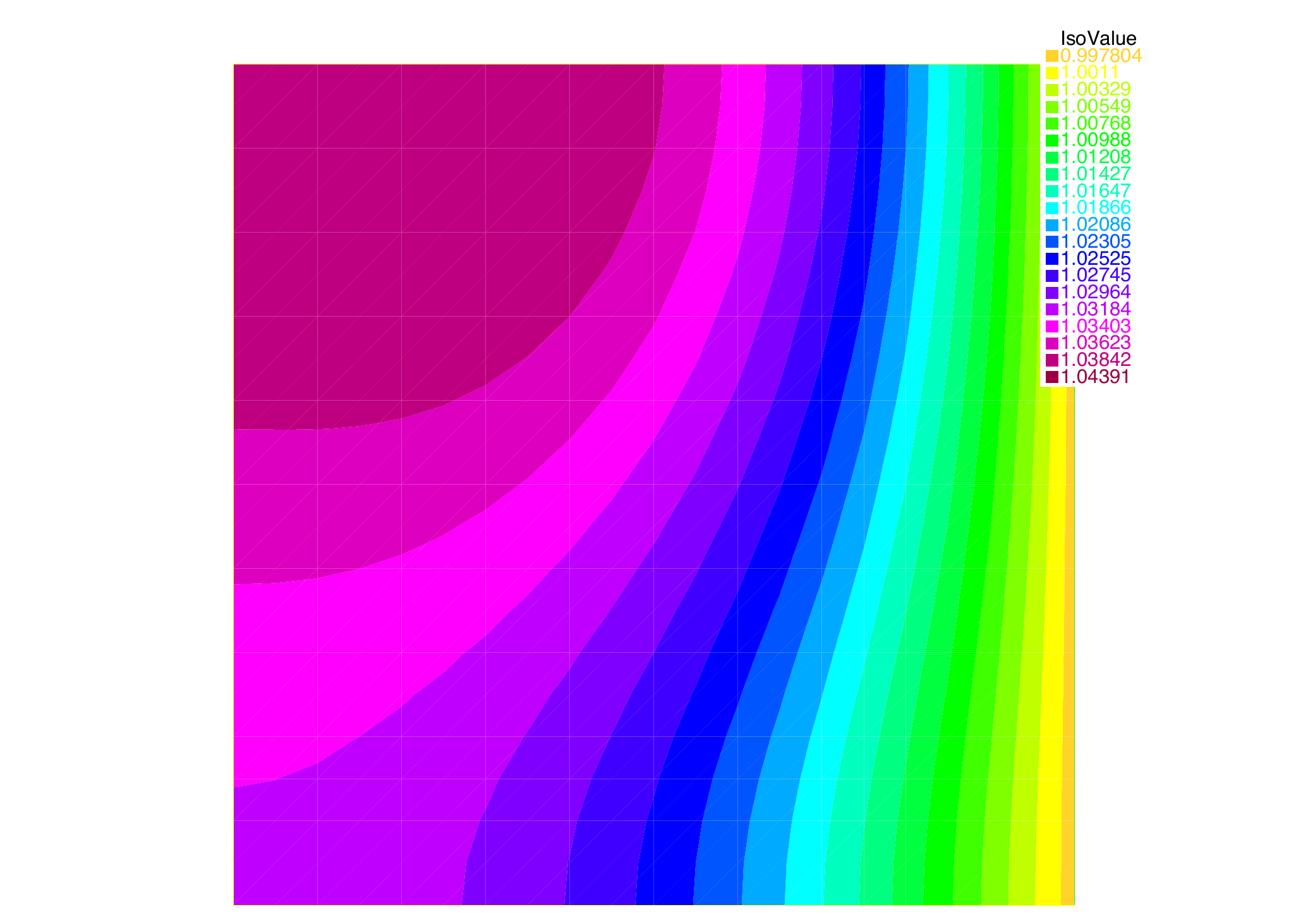}
  \caption{Exact solution of the non-linear Poisson problem (\ref{eq:diffnonlin})}
  \label{Fig-sol2d}
\end{figure}

We now test the non-linear preconditioners on the two dimensional
non-linear diffusion problem (see \cite{fenics})
\begin{equation}\label{eq:diffnonlin}
\left\{\begin{array}{rcll}
-\nabla\cdot ((1+u^2) \nabla u) & = & f,\quad & \Omega=[0,1]^2,\\
u &=& 1,& x=1,\\
\displaystyle\frac{\partial u}{\partial\mathbf{n}} & = &0,&\mbox{otherwise.}
\end{array}\right.
\end{equation}
The isovalues of the exact solution are shown in
Figure \ref{Fig-sol2d}. To calculate this solution, we use a discretization with P1 finite
elements on a uniform  triangular mesh. All calculations have been performed using
FreeFEM++, a C++ based domain-specific language for the numerical solution of PDEs using
finite element methods \cite{freefem}.
We consider a decomposition of the domain into $N\times N$ subdomains
with an overlap of one mesh size $h$, and we keep the number of
degrees of freedom per subdomain fixed in our experiments.  We show in
Table \ref{table:results1lev} a detailed
account of the number of linear
subdomain solves needed for RASPEN and ASPIN at each outer Newton iteration $n$, using the same
notation as in Table \ref{table:results1levFroDisc} (Newton converged in three
iterations for all examples to a tolerance of $10^{-8}$).  We see from these experiments that
RASPEN, which is a non-linear preconditioner based on a convergent
underlying fixed point iteration, clearly outperforms ASPIN, which
would not be convergent as a basic fixed point iteration.

\begin{table}
\centering{
\tabcolsep0.5em
\begin{tabular}{|c|c|c|c|c|c|c|c|c|c|}
\hline
& & \multicolumn{4}{c|}{1-Level} & \multicolumn{4}{c|}{2-Level} \\
\cline{3-10}
$N\times N$ & $n$  &  $ls_n^G$  & $ls_n^{in}$ & $ls_n^{\min}$ &
$LS_n$ & $ls_n^G$  & $ls_n^{in}$ & $ls_n^{\min}$ &
$LS_n$\\\hline
$2\times 2$ & 1 & 15(20) & 4(4) & 3(3) & & 13(23) & 4(4) & 3(3) & \\
                   & 2 & 17(23) & 3(3) & 3(3) & 59(78)& 15(26) & 3(3) & 3(3) & 54(86)\\
                   & 3 & 18(26) & 2(2) & 2(2) & & 17(28) & 2(2) & 2(2) & \\
\hline
$4\times 4$ & 1 & 32(37) & 3(3) & 3(3) & & 18(33) & 3(3) & 3(3) & \\
                   & 2 & 35(41) & 3(3) & 2(2) & 113(132)& 22(39) & 3(3) & 2(2) & 74(126) \\
                   & 3 & 38(46) & 2(2) & 2(2) & & 26(46) & 2(2) & 2(2) &\\
\hline
$8\times 8$ & 1 & 62(71) & 3(3) & 2(2) & & 18(35) & 3(3) & 3(2) & \\
                   & 2 & 67(77) & 3(3) & 2(2) & 211(240) & 23(44) & 3(3) & 2(2) & 77(139)\\
                   & 3 & 74(84) & 2(2) & 1(2) & & 28(53) & 2(2) & 2(1) &\\
\hline
$16\times 16$ & 1 & 125(141) & 3(3) & 2(2) &  & 18(35) & 3(3) & 3(2) &\\
                   & 2 & 136(155) & 2(2) & 2(2) & 418(471)& 23(44) & 2(2) & 2(2) & 75(140)\\
                   & 3 & 150(167) & 2(2) & 1(1) &  & 27(54) & 2(2) & 2(1) &\\
\hline
\end{tabular}
}
\caption{Numerical results with one- and two-level RASPEN and ASPIN for the
  non-linear diffusion problem.}
\label{table:results1lev}
\end{table}

\subsection{A problem with discontinuous coefficients}
\begin{figure}
  \begin{center}
    \includegraphics[width=0.27\textwidth]{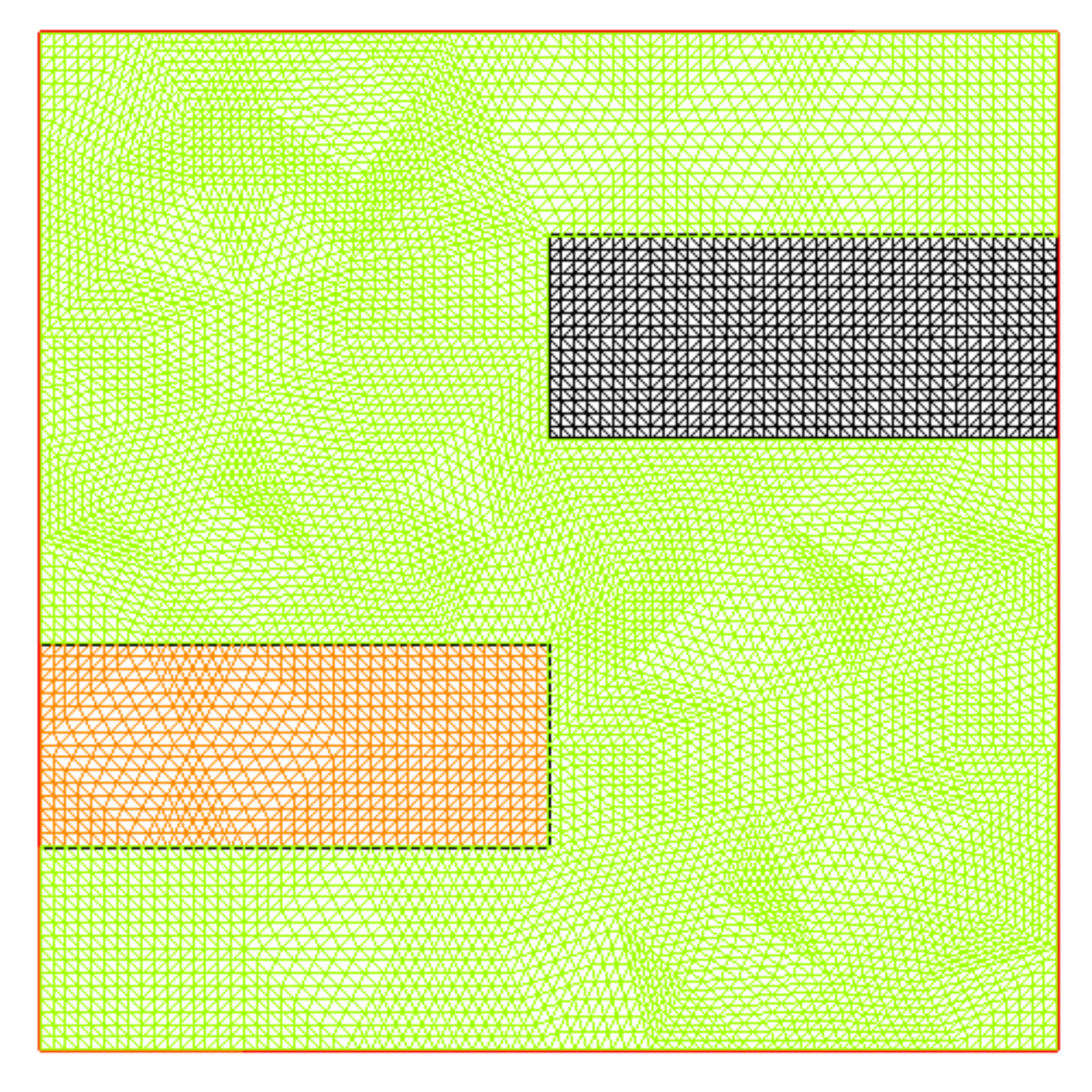}
    \includegraphics[width=0.27\textwidth]{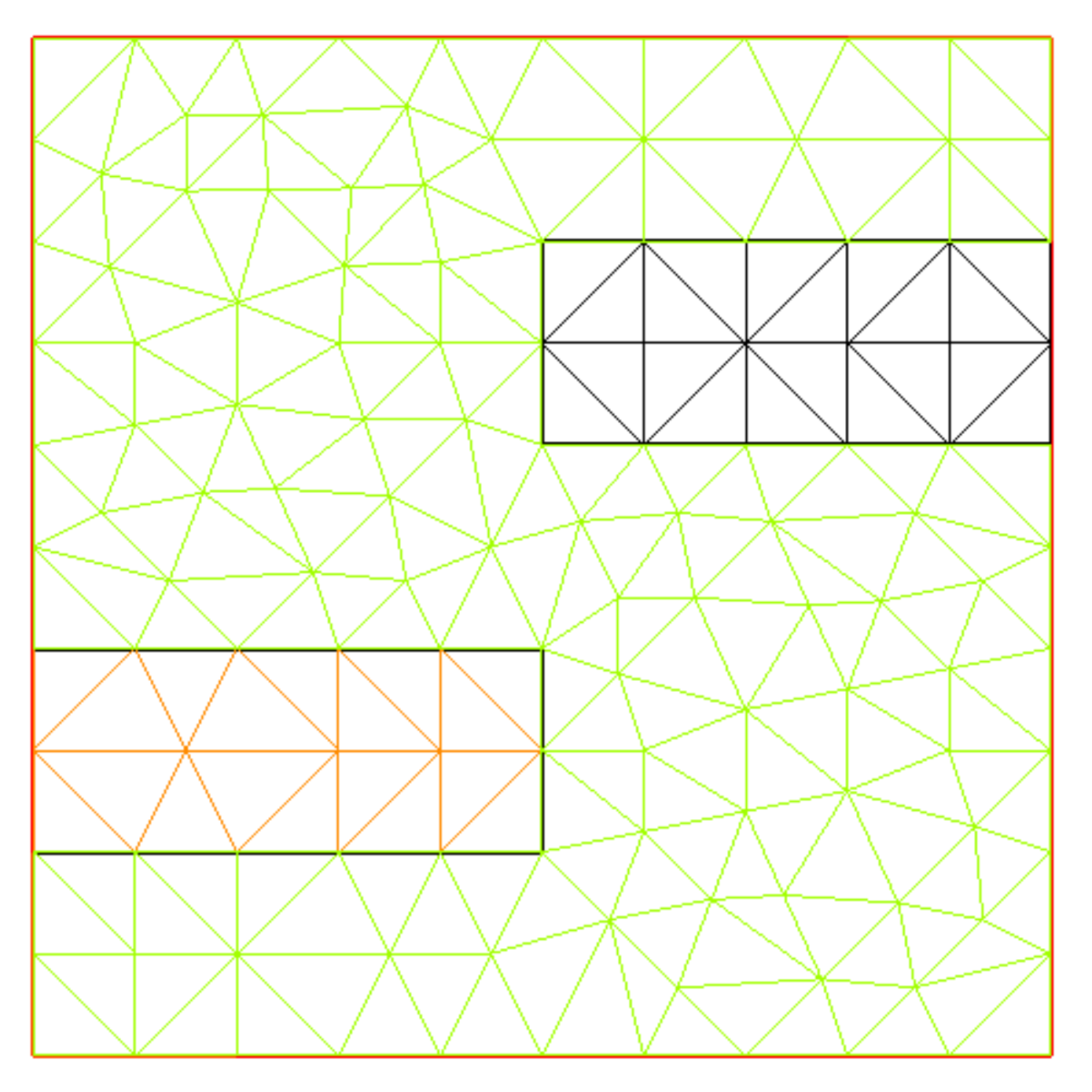}
    \includegraphics[width=0.42\textwidth,bb=6 13 388 253,clip]{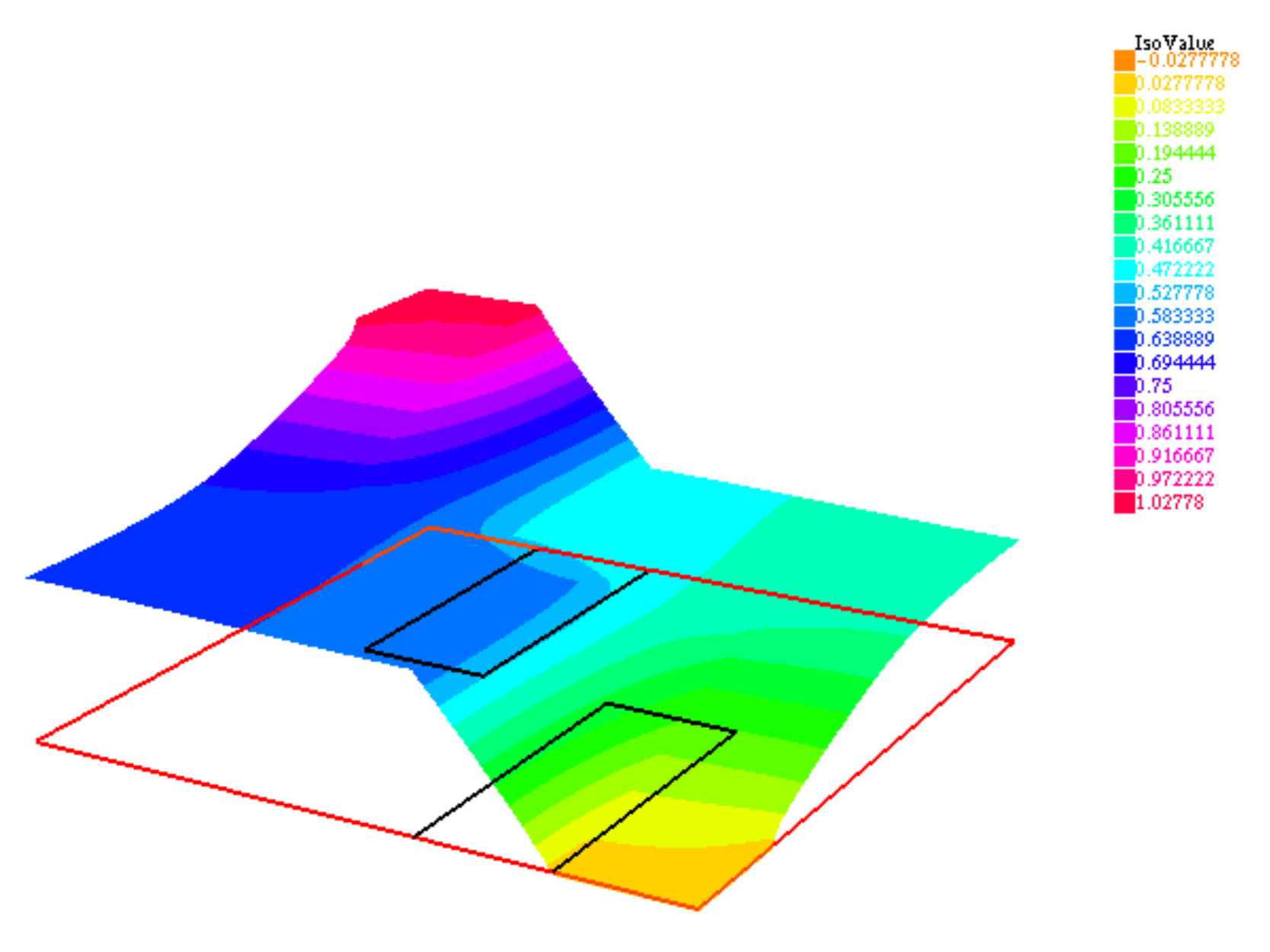}
    \caption{Left: Fine discretization of the domain for $N=4$. The red and black inclusions correspond
      to low-permeability regions. Middle: Coarse grid used for two-level methods.
      Right: Exact solution for the 2D Forchheimer problem.\label{fig:grids}}
  \end{center}
\end{figure}

We now test the non-linear preconditioners on the two dimensional
Forchheimer problem, which can be written as
\begin{equation}\label{eq:2dforch}
\left\{  \begin{gathered}
    -\nabla\cdot {\bf q} = 0,\qquad \Omega=[0,1]^2,\\
    {\bf q} + \beta|{\bf q}|{\bf q} = \Lambda({\bf x})\nabla u,\\
    u = 0\quad \text{on $\Gamma_{d0}$},\qquad u=1\quad \text{on $\Gamma_{d1}$},\\
    {\bf q}\cdot {\bf n} = 0 \quad \text{on $\partial\Omega \setminus (\Gamma_{d0} \cup\gamma_{d1})$,}
  \end{gathered}\right.
\end{equation}
where the Dirichlet boundaries $\Gamma_{d0}$ and $\Gamma_{d1}$ are located at the bottom left and
top right corners of the domain, namely,
$$ \Gamma_{d0} = \{(x,y)\in \partial\Omega : x+y < 0.2\},\qquad \Gamma_{d1} = \{(x,y)\in \partial\Omega: x+y > 1.8\}. $$
The permeability $\Lambda({\bf x})$ is equal to 1000 everywhere except at the two inclusions
shown in red and black in the left panel of Figure \ref{fig:grids}, where it is equal to 1. 
The mesh in the above figure is used to discretize the problem, using P1 finite elements;
the exact solution is shown in the right panel.


\begin{table}
\begin{center}
\begin{tabular}{|l|c|c|c|c|c|c|}
\hline
& \multicolumn{3}{c|}{$\beta=0.1$} & \multicolumn{3}{c|}{$\beta=1$}\\
\cline{2-7}
& $2\times 2$ & $4\times 4$ & $8\times 8$ & $2\times 2$ & $4\times 4$ & $8\times 8$\\
\hline
Newton & 19 & 19 & 19 & 38 & 44 & 48\\
ASPIN & 6 & div. & div.& 6 & div. & div.\\
ASPIN2 & 5 & 6 & 7 & 6 & 7 & 9\\
RASPEN & 5 & 4 & 4 & 5 & 5 & 5\\
RASPEN2 & 4 & 4 & 4 & 5 & 5 & 6\\
\hline
\end{tabular}
\caption{Number of non-linear iterations required for convergence by various algorithms for the 2D Forchheimer problem, as a function of problem size. Divergence of the method is indicated by `div'.\label{nonlineartable}}
\end{center}
\end{table}

\begin{table}
\centering
{\small\tabcolsep0.5em
\begin{tabular}{|c|c|c|c|c|c|c|c|c|c|c|}
  \hline
  & & & \multicolumn{4}{c|}{1-Level} & \multicolumn{4}{c|}{2-Level}\\
  \cline{4-11}
$N\times N$ & $\beta$ & $n$  &  $ls_n^G$  & $ls_n^{in}$ & $ls_n^{\min}$ &
$LS_n$ & $ls_n^G$  & $ls_n^{in}$ & $ls_n^{\min}$ &
$LS_n$ \\\hline
$2\times 2$ & 0.1 & 1 & 22(29) & 5(5) & 5(5) & 82(106) & 10(20) & 6(5) & 5(5) & 47(102) \\
            &     & 2 & 24(32) & 4(4) & 4(4) & & 12(21) & 3(4) & 3(4) &  \\
            &     & 3 & 25(33) & 2(3) & 2(2) & & 14(22) & 2(3) & 2(3) &  \\
            &     & 4 & - (-) &- (-)&- (-)& & - (25) & - (2) &  - (2) & \\
\cline{2-11}
            & 0.2 & 1 & 22(28) & 4(4) & 3(3) & 53(69) & 9(19) & 4(4) & 3(3) & 29(49) \\
            &     & 2 & 24(34) & 3(3) & 3(3) &        & 14(23) & 2(3) & 2(3) &       \\
\cline{2-11}
            & 0.5 & 1 & 22(28) & 4(4) & 4(4) & 53(69) & 9(19) & 4(4) & 3(3) & 29(49) \\
            &     & 2 & 24(34) & 3(3) & 3(3) &        & 14(23) & 2(3) & 2(3) &       \\
\cline{2-11}
            & 1.0 & 1 & 22(28) & 4(4) & 4(4) & 53(69) & 10(21) & 4(4) & 3(3) & 30(51)\\
            &     & 2 & 24(34) & 3(3) & 2(2) &        & 14(23) & 2(3) & 2(2) &       \\
\hline
$4\times 4$ & 0.1 & 1 & 41(53) & 5(5) & 4(4) & 145(179) & 11(21) & 6(6) & 4(4) & 52(111)\\
            &     & 2 & 45(56) & 4(4) & 3(3) &          & 14(23) & 3(4) & 3(3) &        \\
            &     & 3 & 48(58) & 2(3) & 2(2) &          & 16(24) & 2(4) & 2(3) &        \\
            &     & 4 & - (-) & - (-)&- (-)&          & - (26) & - (3) & - (2) & \\
\cline{2-11}
            & 0.2 & 1 & 41(52) & 4(4) & 3(3) & 94(118) & 11(21) & 4(4) & 3(3) & 33(54) \\
            &     & 2 & 47(59) & 2(3) & 2(2) &         & 16(26) & 2(3) & 2(2) &        \\
\cline{2-11}
            & 0.5 & 1 & 41(51) & 4(4) & 3(3) & 94(116) & 11(21) & 4(4) & 3(3) & 33(54) \\
            &     & 2 & 47(58) & 2(3) & 2(2) &         & 16(26) & 2(3) & 2(2) &        \\
\cline{2-11}
            & 1.0 & 1 & 41(51) & 4(4) & 3(3) & 94(116) & 11(21) & 4(3) & 3(3) & 34(53) \\
            &     & 2 & 47(58) & 2(3) & 2(2) &         & 17(26) & 2(3) & 2(2) &        \\
\hline
$8\times 8$ & 0.1 & 1 & 86(104) & 5(5) & 3(3) & 468(573) & 16(24) & 6(5) & 3(3) & 73(160) \\
            &     & 2 & 92(111) & 3(4) & 2(2) &          & 21(27) & 4(4) & 3(3) &       \\
            &     & 3 & 95(115) & 3(3) & 2(2) &          & 24(26) & 2(4) & 2(2) &       \\
            &     & 4 & 90(116) & 2(2) & 1(1) &          & - (30) & - (3) & - (2) &     \\
            &     & 5 & 90(111) & 2(2) & 1(1) &          & - (35) & - (2) & - (1) &     \\
\cline{2-11}
            & 0.2 & 1 & 84(103) & 4(4) & 3(3) & 373(457) & 16(24) & 4(4) & 3(3) & 46(62)\\
            &     & 2 & 93(115) & 3(3) & 2(2) &          & 24(31) & 2(3) & 2(2) &       \\
            &     & 3 & 94(117) & 2(2) & 1(1) &          & - (-) &- (-)&- (-)&       \\
            &     & 4 & 91(111) & 2(2) & 1(1) &          & - (-) &- (-)&- (-)&       \\

\cline{2-11}
            & 0.5 & 1 & 84(104) & 4(4) & 3(3) & 374(461) & 16(25) & 4(4) & 3(3) & 46(63) \\
            &     & 2 & 94(115) & 3(3) & 2(2) &          & 24(31) & 2(3) & 2(2) &       \\
            &     & 3 & 94(119) & 2(2) & 1(1) &          & - (-) &- (-)&- (-)&      \\
            &     & 4 & 91(112) & 2(2) & 1(1) &          & - (-) &- (-)&- (-)&      \\
\cline{2-11}
            & 1.0 & 1 & 84(104) & 4(4) & 2(2) & 375(461) & 16(25) & 4(4) & 3(2) & 47(64) \\
            &     & 2 & 95(115) & 3(3) & 2(2) &          & 25(32) & 2(3) & 2(2) &       \\
            &     & 3 & 95(119) & 2(2) & 1(1) &          & - (-) &- (-)&- (-)&      \\
            &     & 4 & 91(112) & 2(2) & 1(1) &          & - (-) &- (-)&- (-)&      \\
\hline
\end{tabular}
}
\caption{Numerical results with one- and two-level RASPEN and ASPIN for the
  2D Forchheimer problem. `-' indicates that the method has converged.}
\label{table:results2dforch1lev}
\end{table}

We consider a decomposition of the domain into $N\times N$ subdomains
with an overlap of one mesh size $h$, and we keep the number of
degrees of freedom per subdomain fixed in our experiments.  For the two-level methods,
the coarse function $F_0$ consists of a P1 discretization
of the problem over the coarse grid shown in the middle panel of Figure \ref{fig:grids}.

\FK{To measure the difficulty of this problem, we run our non-linear algorithms (standard
Newton, one and two-level ASPIN, one and two-level RASPEN) on this problem for $\beta=0.1$
and $\beta=1$. We show in Table \ref{nonlineartable} the number of iterations required
for the convergence of each algorithm. We see that between the discontinuous permeability
and the non-linearity introduced by $\beta$, standard Newton requires many iterations to
converge, and one-level ASPIN diverges for the larger problems. On the other hand, one and
two-level RASPEN (and two-level ASPIN, to a lesser extent) converge in a small number of
non-linear iterations.}

\FK{Next, we compare the one and two-level variants of ASPIN and RASPEN in terms of the
total amount of computational work. To deal with the convergence problem in one-level
ASPIN, we adopt the continuation approach,}
where we solve the
problem for a sequence of $\beta$ (0, 0.1, 0.2, 0.5 and 1.0), using the solution for the
previous $\beta$ as the initial guess for the next one.
Table \ref{table:results2dforch1lev} shows a detailed
account for each outer Newton iteration $n$ of the linear
subdomain solves needed for both RASPEN and ASPIN using the same
notation as in Table \ref{table:results1levFroDisc}. We omit the data for $\beta = 0$, as the problem
becomes linear in that case.
\FK{We see again from these experiments that the
RASPEN-based preconditioners can handle non-linearly difficult problems,
requiring fewer non-linear iterations and linear solves than their ASPIN counterparts.}

\section{Conclusion}

We have shown that just as one can accelerate stationary iterative
  methods for linear systems using a Krylov method, one can also
  accelerate fixed point iterations for non-linear problems using
  Newton's method. This leads to a guiding principle for constructing
  non-linear preconditioners, which we illustrated with the systematic
  construction of RASPEN. While this design principle leads to good
  non-linear (and linear) preconditioners, \MG{see for example
  \cite{haeberlein2013newton,haeberlein2015schwarz} for a similar approach for non-linear
  evolution problems}, \FK{it is by no means the only approach possible;
  in the linear case, for instance, the additive Schwarz preconditioner \cite{Dryja:1987:AVS}, 
  as well as
  the highly effective and robust FETI preconditioner \cite{feti}
  and its variants, do not correspond to a convergent iteration. Indeed,}
  clustering the spectrum into a few clusters is sometimes better than having a small
  spectral radius, see for example the results for the HSS preconditioner in
  \cite{benzi2003optimization}. Thus, it is still an open question whether there are other
  properties that a preconditioner should have that would make it more
  effective, \FK{even if it is associated with a divergent iteration}.  For non-linear preconditioning, maybe
  it is possible to greatly increase the basin of attraction of the
  non-linearly preconditioned Newton method, or to improve its
  preasymptotic convergence, before quadratic convergence sets in. It
  also remains to carefully compare RASPEN with linear preconditioning
  inside Newton's method; promising results for ASPIN can be found
  already in \cite{skogestad2013domain}.

\section*{Acknowledgments} We would like to thank TOTAL for partially supporting this work.

\bibliographystyle{abbrv}
\bibliography{paper_revised}

\end{document}